\documentclass[reqno,11pt]{amsart}
\usepackage{amsmath, amssymb, amsthm,amsfonts} 
\usepackage[english]{babel}
\usepackage{bbm,bm}
\usepackage{graphicx}
\usepackage{url}
\usepackage{epstopdf}
\usepackage[a4paper,margin=1in]{geometry}
\usepackage{subcaption}
\usepackage{textcomp}
\usepackage[ruled,vlined]{algorithm2e}
\usepackage{multicol}
\usepackage{multirow}
\usepackage{wrapfig}
\usepackage{diagbox}
\usepackage{floatrow}
\usepackage{mathtools}
\usepackage[comma]{natbib}
\usepackage[linguistics]{forest} 
\usepackage[titletoc,toc,title]{appendix}
\usepackage{setspace}
\newfloatcommand{capbtabbox}{table}[][\FBwidth]
\usepackage{abstract}
\usepackage{tikz}
\usetikzlibrary{arrows, automata,positioning,calc,shapes,decorations.pathreplacing,decorations.markings,shapes.misc,petri,topaths,matrix}
\usepackage{pgfplots}
  \pgfplotsset{compat=newest}
  \usetikzlibrary{plotmarks}
  \usepackage{grffile}
\newlength\figureheight
  \newlength\figurewidth
\pgfplotsset{%
    tick label style={font=\scriptsize},
    label style={font=\footnotesize},
    legend style={font=\footnotesize},
         every axis plot/.append style={very thick}
}

\usepackage{rotating}
\usepackage{amsbsy,enumerate}
\usepackage{graphicx}
\usepackage{comment}
\usepackage{mathrsfs} 
\usepackage{tikzscale}
\usepackage{xcolor}

\newcommand{\vb}{\vspace{5mm}}

\newcommand\addtag{\refstepcounter{equation}\tag{\theequation}}

\allowdisplaybreaks

\makeatletter
\newtheorem*{rep@theorem}{\rep@title}
\newcommand{\newreptheorem}[2]{%
\newenvironment{rep#1}[1]{%
 \def\rep@title{#2 \ref{##1}}%
 \begin{rep@theorem}}%
 {\end{rep@theorem}}}
\makeatother

\makeatletter
\renewcommand*\env@matrix[1][\arraystretch]{%
  \edef\arraystretch{#1}%
  \hskip -\arraycolsep
  \let\@ifnextchar\new@ifnextchar
  \array{*\c@MaxMatrixCols c}}
\makeatother

\newreptheorem{theorem}{Theorem}

\newtheorem{remark}{Remark}
\newtheorem{example}{Example}

\newreptheorem{proposition}{Proposition}

\begin{document}

\doublespacing

\title[Optimal Departure Time Advice]{Optimal Departure Time Advice in Road Networks \\ with Stochastic Disruptions}
\author{Rens Kamphuis, Nikki Levering and Michel Mandjes}

\maketitle

\begin{abstract}
Due to recurrent (e.g.\ daily or weekly) patterns
and non-recurrent disruptions (e.g.\ caused by incidents), travel times in road networks 
are time-dependent and inherently random.
This is challenging for travelers planning a future trip,
aiming to ensure on-time arrival at the destination,
while also trying to limit the total travel-time budget spent.
The focus of this paper lies on determining 
their \textit{optimal departure time}:
the latest time of departure for which a chosen 
on-time arrival probability can be guaranteed.
To model the uncertainties in the network,
a Markovian background process is used, tracking
events affecting the driveable vehicle speeds on
the links, thus {enabling} us to incorporate both recurrent and non-recurrent effects.
It allows the evaluation of the travel-time distribution,
given the state of this process at departure,
on each single link.
% Then a computationally efficient algorithm is devised that uses these individual 
% link travel-time distributions to obtain the optimal departure time and corresponding
% optimal path for a trip from the origin to a given destination.
{Then, a computationally efficient algorithm is devised that uses these individual 
link travel-time distributions to obtain the optimal departure time for a given path
or origin-destination pair.}
Since the conditions in the road network, and thus
the state of the background process, may change between the time of request 
and the advised time of departure, we consider
an online version of this procedure as well, in which the traveler
receives departure time updates while still at the origin. 
Finally, numerical experiments are conducted
to exemplify a selection of properties of the optimal departure time and, moreover, quantify the performance of the presented algorithms in an existing road network --- the Dutch highway network.

% Now, based on the state of this process at the moment of request,
% the resulting travel time distribution for a given path
% and departure time can be computed via
% a discretization procedure.
% First, we present a discretization procedure to 
% obtain the resulting travel time distribution for
% a given departure time and path.
% Then, we show how to use a bisection procedure to
% obtain the optimal departure time for a given OD-pair.
% To this end, we present both an exact but computationally
% demanding A$^\star$-algorithm, and the $k$-shortest path.

% Due to recurrent daily patterns
% and non-recurrent disruptions, travel times in road networks 
% suffer from both time-dependence and randomness.
% We focus on a departure time problem: our objective
% is to find, for a single vehicle planning a trip from
% an origin to a destination, the latest departure time for which a chosen
% on-time arrival probability can be guaranteed.
% This is a challenge often faced by travelers planning a future trip,
% aiming to ensure on-time arrival,
% whilst also trying to limit the total travel time budget.
% To model the uncertainties in the network,
% we use a Markovian background process that tracks
% recurrent and non-recurrent events affecting
% the speed of travelers. We present a discretization
% procedure to obtain travel time distributions.

% To model these disruptions, we use a Markovian 
% background process
% that tracks the occurrence of incidents
% and their effect on the arc speeds.
% The resulting velocity dynamics describe a travel
% time distribution for the traveler planning the trip.

\vb


\vb

\noindent
{\sc Affiliations.} 
{RK, NL and MM are with the Korteweg-de Vries Institute for Mathematics, University of Amsterdam, Amsterdam. MM is also with E{\sc urandom}, Eindhoven, and with Amsterdam Business School, University of Amsterdam, Amsterdam.}
This research project is partly funded by the NWO Gravitation project N{\sc etworks}, grant number 024.002.003. Date: \today. 

\end{abstract}

\newpage

\section{Introduction} \label{sec: introduction}

{\sc Motivation} 

Travel times in road networks are inherently random. 
%A traveler who wants to reach a given destination from their current location, does not know the arrival time with certainty. 
This  could lead to late arrival at the destination, which, depending on the importance of the trip, may be highly undesirable. At the same time, travelers will generally also be reluctant to depart far too early, as this would leave them with less time available for other  activities. Thus, there is a clear need for algorithms that generate the latest departure time such that a certain user-specified on-time arrival probability can be guaranteed. 
% \textcolor{red}{This demand
% exists both for travelers that have specified the path to travel, 
% as well as for users that have not selected a route, but only an origin-destination (OD) pair.
% }
% \textcolor{red}{This demand
% exists for travelers that have specified the path to travel, 
% as well as for users that have only selected their destination
% and are free to travel via any path between the
% origin-destination (OD) pair.
% }

The randomness of travel times is due to 
%the level of congestion being influenced by 
both recurrent patterns (e.g., peak hours) and non-recurrent events (e.g., traffic accidents, vehicle failures, unfavorable weather conditions). Whereas recurrent patterns are predictable from historic data, the time and location of future non-recurrent events are intrinsically uncertain. It is noted, though, that the {\it current} state of the road network {\em is} known; think of the locations of the incidents that are currently present, as well as their elapsed durations. As currently present incidents further away from the origin are more likely to be cleared once the corresponding road segment is reached by the traveler, this information can be exploited when determining the optimal departure time.

When devising a procedure to generate an optimal departure time {for a given path or origin-destination (OD) pair}, 
various considerations play a role. 
Ideally, one would like to determine a departure time that (i)~takes the risk averseness of the traveler into account, (ii)~incorporates the randomness of travel times, due to both recurrent and non-recurrent effects, and (iii)~exploits knowledge of the locations of the currently present non-recurrent events. To be able to discuss our contributions, we proceed by an account of the existing literature (where we do not attempt to provide an exhaustive overview, but primarily focus on contributions relevant in the context of our work).

{\sc Literature} 

While there are different notions of optimality, optimal reliable paths can roughly be classified into two categories. Either a path is regarded as optimal if it has the highest on-time arrival probability, as in the seminal work by \cite{frank1969shortest}, or, a path is optimal if it has the lowest expected travel time while also guaranteeing a certain reliability threshold. An example of the latter is the $\alpha$-shortest path as introduced in \cite{ji2005models}, which is the path with the lowest expected travel time such that the probability of arriving on time is at least $\alpha$. A study to identifying the a-priori $\alpha$-shortest path can be found in \cite{Nie2009ShortestProbability}, whereas \cite{fan2005arriving} and \cite{christman2013maximizing} outline a dynamic routing procedure to maximize the on-time arrival probability.

% The earliest work on routing in graphs with stochastic edge lengths can be found in \cite{frank1969shortest}. This work is concerned with finding the probability distribution of the shortest path, which is then later used for a pairwise path comparison test. Here, optimality entails the path that maximizes the probability of realizing a weight less than a given threshold. A different notion of optimality is used in \cite{sigal1980stochastic}, where a path is regarded as optimal if it has the greatest probability of realizing the least weight.

In \cite{loui1983optimal} it is recognized that drivers may have different preferences, and that the optimal route is the route that optimizes the personal utility function of the individual driver. Often, these utility functions  take into account both the expected travel time and the reliability of the route; see e.g.\ \cite{sen2001mean}. The typical objective is to identify the route that mimimizes the sum of the expected travel time and a term that quantifies the travel-time variance, with weights that reflect the driver's risk averseness. Such a utility function is intended to strike a proper balance between minimizing the  travel time and controlling its randomness  \citep{tilahun2010moment}.
The previously mentioned $\alpha$-shortest path criterion succeeds in favoring routes with a low mean travel time, while also penalizing travel-time uncertainty. As this criterion uses the full travel-time distribution (rather than just the travel time variance), the uncertainty of a route can be captured in more detail than in the approach of \cite{sen2001mean}. 
Moreover, the degree of the uncertainty penalization can be adjusted by the parameter $\alpha$, which is chosen to reflect the specific driver's risk averseness.

An important observation is that the literature on optimal reliable routing,
including the above-mentioned works, typically assumes instant departure of the traveler. 
However, besides late arrival, also early arrival evidently incurs an (opportunity) cost, so that immediate departure is generally not optimal \citep{small1982scheduling}. Whereas maximizing the on-time arrival probability would require instant departure, the $\alpha$-shortest path criterion can still be used in case of non-instant departure. The idea is then to identify the latest departure time such that an on-time arrival probability of at least $\alpha$ can be guaranteed. This approach is followed in e.g.\  \cite{chen2014reliable}, where besides an optimal route, also the optimal departure time is determined under a stochastic first-in-first-out assumption; see also  \cite{yang2020finding}, assuming
log-normally distributed travel times.

As mentioned, travel times cannot be assumed constant; sources of these fluctuations are discussed in \cite{stogios2014incorporating}. To deal with the recurrent patterns, one could work with the framework of \cite{chen2014reliable} and \cite{yang2020finding}, where time is divided into intervals and the parameters of the travel-time distributions have interval-specific values. These parameters 
have been inferred from historic data, and hence are considered known. However, travel times are also subject to non-recurrent effects, such as congestion due to incidents. These should clearly be taken into account, as the occurrence or clearance of incidents greatly affects the travel-time distribution. Obviously, one cannot use techniques designed to deal with recurrent patterns, but one would like to exploit knowledge of the statistics of incident durations and inter-incident times, as well as knowledge of the current state of the road network (in particular the locations of incidents that are currently present).

{\sc Contributions} 

It is the main aim of this paper to contribute to the literature on determining optimal departure times, with a focus on including the effect of non-recurrent congestion, aiming at meeting the three requirements (i), (ii), and (iii) identified above. A more detailed account of the contributions of our work is as follows:

\begin{itemize}
   \item[$\circ$]
To model the travel-time dynamics such that it can cover both recurrent and non-recurrent effects, we use the Markovian Velocity Model ({\sc mvm}) as introduced in \cite{Levering2021AConditions}. 
The 
{`velocity oriented' and data-driven}
{\sc mvm}, building upon the stochastic disruption models in 
\cite{Psaraftis1993,Karoufeh2004,Kim2005a,Kim2005b,Thomas2007,Guner2012}, and \cite{Sever2013}, uses a Markovian background process that tracks the events that affect the driveable vehicle speeds on road segments.
% \textcolor{red}{
% This framework, which is `velocity oriented' and data-driven, offers extreme flexibility for the modeling of travel time distributions.
% In addition,} 
A particularly attractive aspect of the {\sc mvm} is that knowledge of the current state of the road network (e.g.\ locations of incidents) can be naturally incorporated. %Thus, besides capturing non-recurrent events, it also succeeds in incorporating readily available information about the current state of the network.

\item[$\circ$]
The velocities in the road network are fully determined by the current state of the Markovian background process. The dynamics of the background process determine the travel-time distribution for each link, which can then be used to obtain the travel-time distribution for any route in the network. Having access to the distribution, %besides computing the expected travel time, 
the probability of on-time arrival for any departure time can be computed.
%as well. 
This allows us to set up a procedure that identifies, for any desired on-time arrival probability, the optimal departure time.
{We do this both for a specified path and for 
a specified origin-destination pair (i.e., also providing
the path that yields the optimal departure time for this OD pair).}
% \textcolor{red}{We consider
% both travelers that have specified the path to travel, 
% as well as users that have only selected their destination
% and are free to travel via any path between the
% origin-destination (OD) pair.
% }
The risk averseness is naturally taken into account, as the departure time depends on the selected value of the on-time arrival probability.

\item[$\circ$]
Where \cite{chen2014reliable} and \cite{yang2020finding} exclusively rely on historic estimates to evaluate the travel-time distribution, our procedure also includes the impact of the current state of the road network. Rather than choosing a specific travel-time distribution, such as the log-normal one featuring in  \cite{yang2020finding}, the {\sc mvm} is highly flexible, as extensively argued in \cite{Levering2022Estimation}, in particular in terms of modeling the distribution of the incident durations and the inter-incident times. Another advantage of our approach is that the road network by construction satisfies the desirable FIFO property, and we do not have to rely on a stochastic-FIFO assumption as in \cite{chen2014reliable}.

\item[$\circ$]
We recognize that the conditions in the road network, and thus the state of the Markovian background process, may change between the time of request and the advised time of departure. Therefore, we also consider an online version of the optimal departure time problem, in which the traveler receives departure time updates while still at the origin.

\item[$\circ$]
{Lastly, we provide a selection of numerical experiments in which we implement our procedure and study the properties of the optimal departure time. We show that the optimal departure time is greatly affected by both the state of road network and the time of request. Moreover, we quantify the (potentially substantial) gains in travel time budget that can be realized by utilizing the online version of the optimal departure time problem. We also examine the efficiency of our procedure, and demonstrate that our procedure can successfully be employed in a real-world road network.}

\end{itemize}

This paper is organized as follows. Section~\ref{sec:problem} describes our model and objective. Then, in Section~\ref{sec:optdeptime} a procedure 
is presented for obtaining the travel-time distribution of individual
links, after
which computationally efficient algorithms are outlined that
output the optimal departure time for a given path {or OD pair}, target arrival time and on-time probability. 
Numerical experiments are performed in Section~\ref{sec:numex}, that aim to exemplify a selection of properties of the optimal departure time and studies the efficiency of our procedure. Finally, Section~\ref{sec:concl} gives some conclusive remarks.
{
Various extensions to our procedure are discussed in the appendix.}

\section{Model and problem description} \label{sec:problem}

% Introductory paragraphs -> potentially to introduction
We consider a single vehicle 
that plans an upcoming trip between 
an OD-pair
in a road network, wishing to arrive
at the destination before a given time.
Travel times in the network are subject to
stochastic disruptions, i.e., events that
affect driveable vehicle speeds.
These disruptions (in the sequel often referred to as incidents)
are modeled via a Markovian 
background process.
% and their effect on the arc speeds.
% The resulting velocity dynamics describe travel
% time distributions for vehicles in a
% road network.
At the time the traveller requests advice (with which we identify $t=0$ throughout this paper), the state
of the background process is known. If one would have access to the travel-time distribution for any departure time (after $t=0$, that is), given knowledge of the background process at $t=0$, one could find the \textit{optimal departure time}, i.e., the latest departure time such that the probability of arriving on-time is at least, say, $\eta$. A subtlety is that it should be incorporated that the state of the background process could change between the request time and the departure time.

% The resulting velocity dynamics describe a travel
% time distribution for the traveler planning the trip.
% Now, knowing this travel time distribution, 

In Section~\ref{subsec: Markov} we describe the Markovian
Velocity Model ({\sc mvm}) employed to model the effect of
disruptions on the vehicle
speeds. 
With the notation introduced there,
Section~\ref{subsec: departureproblem} formally introduces the
problem of determining the optimal departure time.

\subsection{Travel Time Dynamics: Markovian Velocity Model} \label{subsec: Markov}~\\
Let $G = (N,A)$ be a graph representation of the 
considered road network, of which the set of
nodes $N$ represents the intersections in the road 
network and the set of directed arcs $A$ represents
the roads connecting these intersections. Thus,
we have that $(k,\ell) \in A$ if and only if there
is a directed road segment from the intersection
represented by node $k$ to the intersection represented
by node $\ell$. Throughout, we let $n := |A|$ and 
write $A = \{a_1,\dots,a_n\}$ for the set of arcs in $G$,
with $a_i := (k_i,\ell_i)$ for some $k_i,\ell_i \in N$.

In reality, when traversing a link $a \in A$,
the driveable vehicle speed is not always constant,
as the occurrence of traffic events potentially affects
the velocities.
Importantly, the most prominent source of speed variability
is formed by random traffic incidents, such as
accidents, on-road obstacles, and vehicle break-downs.
The {\sc mvm}, originally
introduced in \cite{Levering2021AConditions}, uses
a Markovian background process $B(t)$ to record the 
evolution of these traffic incidents in the road network.
{By doing so, the model provides a direct relationship between traffic events and travel times,
creating a transparent modelling framework that can be made operational with relatively low effort (as demonstrated in \cite{Levering2022Estimation}).}
% As stated in Section~\ref{sec: introduction}, 
% we originally introduced the MVM in \cite{Levering2021AConditions}, 
% in which the resulting
% travel time dynamics are used 
% to route with
% minimum expected travel time.

For $a_i \in A$, define $\{X_{a_i}(t), t \geqslant 0\}$
as a two-state continuous Markov process with
\begin{align} 
    X_{a_i}(t) &= 
    \begin{cases}
    1 &\quad \quad \text{if there is no incident on arc $a_i$ at time $t$} \\
    2 &\quad \quad \text{otherwise}.
    \label{eq:twostate}
    \end{cases}
\end{align}
Hence, $X_{a_i}$ provides information about a possible
incident on the arc $a_i$ in $G$.
We let the velocity of a vehicle
traversing $a_i$ depend on the occurrence of incidents
in the following way: the vehicle speed at time $t$ 
equals $v_{a_i}(s_i)$ if $X_{a_i}(t)$ is in state $s_i \in \{1,2\}$.
We assume that for $i$ unequal to $j$ 
the processes $X_{a_i}(t)$ and $X_{a_j}(t)$
evolve independently. Consequently,  
for $i = 1,\dots,n$, the process $X_{a_i}(t)$ 
is completely described by its transition rate matrix
\begin{equation}
    Q_{a_i} = 
    \begin{bmatrix}
    -\alpha_i & \alpha_i \\ \beta_i & -\beta_i
    \end{bmatrix},
    \label{eq:qmatrix}
\end{equation}
with $\alpha_i,\beta_i \in \mathbb{R}_{>0}$,
and initial state $X_{a_i}(0)$.
Thus, observe that $X_{a_i}(t)$ is a process
that cyclically switches between an incident
and incident-free state, whose durations
have exponential distributions with
mean $1/\alpha_i$ and $1/\beta_i$ respectively.
In conclusion,  
\begin{align*}
    B(t) := (X_{a_1}(t),\dots,X_{a_n}(t)),
\end{align*}
is a Markovian background process with state space
$S = \{1,2\}^n$, that tracks the occurrence of incidents
in the road network. The {\sc mvm} then describes 
the effect of this background process, characterized
by its initial $2^n$-dimensional state $B(0)$ and $(2^n\times 2^n)$ transition rate matrix denoted $Q$,
on the arc speeds.

\begin{remark}{\em 
It is noted that the {\sc mvm} as presented in this subsection does not
cover the full potential of the {\sc mvm} as demonstrated
in \cite{Levering2021AConditions}.
However, as the main focus of this paper regards the deduction
of the optimal departure time in a network
with stochastic velocities, and not the model
of these stochastic velocities itself, 
we use, for expositional reasons, a compact
version of the {\sc mvm}.
{Specifically, in this paper we 
only allow the speed on arc $a_i$
to depend on the state of $X_{a_i}(t)$, 
whereas the full
{\sc mvm} allows the speed to depend on the complete state of $B(t)$.
Moreover,  
we only consider the processes $X_{a_i}(t)$
in the two-state form~\eqref{eq:twostate}, whereas
the full {\sc mvm} allows any continuous-type Markov process.
Thus, by working with phase-type distributions, which are capable of approximating
any distribution (on the positive numbers) arbitrary closely, also non-exponential incident durations and inter-incident
times could be incorporated.
Lastly, whereas we only consider the impact of incidents,
the process $B(t)$ can easily be extended with a
Markov process $Y(t)$ that captures the recurrent traffic patterns.
Importantly, the presented results and algorithms
can be extended to include these three generalizations; for the latter two generalizations, this is described in more detail in Appendix~\ref{app:recurrent}}. 
% Also non-exponential incident durations and inter-incident times can be incorporated; this is done by working with phase-type distributions, which are capable of approximating any distribution arbitrarily closely.
}\hfill$\Diamond$
\end{remark}

The {\sc mvm}, 
particularly in its full version (Remark~1), is
extremely flexible, and therefore well capable of
describing real-world travel-time distributions.
Indeed, in \cite{Levering2022Estimation} it is shown
how a database of incident registrations 
and loop detector data can be used
to operationalize the {\sc mvm}. Specifically,
it is shown how to model the 
randomness of incidents lengths
and inter-incident durations, and set
the corresponding driveable speed levels.
This allows the fitting of the {\sc mvm}
in e.g.\ the Dutch road network, in which 
a high density of loop detectors and 
lists of occurred incidents are available.
% Now, before being able to use the {\sc mvm} for obtaining
% the optimal departure time, such a fitting
% procedure should be executed.
% The operationalized {\sc mvm} can then be
% used for various applications, such as
% routing with minimum expected travel time
% \citep{Levering2021AConditions}, and, as we
% argue 

\subsection{Objective: Departure Time Advice} \label{subsec: departureproblem}~\\
In this subsection we outline the problem of determining the optimal departure time. In doing so, we distinguish
between two settings, called the {\it offline} and the {\it online} setting. In the offline setting, the optimal departure time is only determined once, based on the information available at the request time, after which the traveler will indeed leave at this time instance. In the online setting, however, the departure time can be updated while the traveler is waiting, as changing conditions may result in a new optimum. 

\medskip

\textit{Offline Setting --- }
We will first formally define the optimal departure time. Suppose a traveler requests a route at the current time, i.e, $t = 0$, and is interested in the latest time of departure such that a certain arrival on-time probability can be guaranteed. If the requested arrival time is given by $t = M>0$, i.e. $M$ time units from the time of request, and if the desired on-time arrival probability is at least $\eta\in(0,1)$, the optimal departure time of a traveler is defined as
\begin{equation}
	t^\star  := \sup\big\{t \geqslant  0 :\, \mathbb{P}\big(Y_{t} \leqslant M \mid B(0)\big) \geqslant  \eta\big\}.
	\label{optimal_departure_offline}
\end{equation}
In \eqref{optimal_departure_offline}, the random variable $Y_{t}$ represents the arrival time  when departing at time $t$, and can be written as $Y_{t} = t + T_{t}$, with $T_{t}$ the travel time when departing at time $t$. Importantly, 
the random variables $Y_{t}$ and $T_{t}$ are affected by the current background state $B(0)$. It is noted that the distributions of $Y_{t}$ and $T_{t}$ will change once future information about the state of the network becomes available (i.e., $B(s)$, when $s$ is approaching $t$). In this offline setting, however, we  assume that only the current state of the network $B(0)$ is known. Note that the on-time arrival probability at time of departure will generally differ from $\eta$.
Of course, it may happen that there exists no $t\geqslant  0$ that satisfies the condition in \eqref{optimal_departure_offline}. In this case, we put $t^\star  := -\infty$, and it depends on the preferences of the driver to either depart immediately or to not depart at all.

We will now pay closer attention to the conditional probability in \eqref{optimal_departure_offline}. 
First, recognize that the travel-time distribution depends on 
the departure time $t\geqslant 0$
only through the state $B(t)$, which is unknown at time 0.
%At $t = 0$, the state of $B(t')$ at departure time $t' > 0$ is unknown. 
However, it is possible to determine the distribution of $B(t)$ by using the known current state $B(0) = s$ and transition matrix $Q$. 
{Using general results for continuous-time Markov chains
\citep{Norris1997, Ibe2013}}, it follows that this distribution is given by
\begin{equation}
{\boldsymbol p}_{t}^s := \left(\,\mathbb{P}(B(t) = s' \,|\, B(0) = s)\,\right)_{s' \in S} = {\boldsymbol p}_{0}^s\, e^{Q t},
\label{distBOffline}
\end{equation}
with ${\boldsymbol p}_{0}^s$, by definition, an row vector of dimension $|S|=2^n$, with a 1 at the entry that corresponds to the state $s \in S$ and a 0 at every other entry. Thus, with $T[{s'}]$ denoting the travel time corresponding to departing when the background
process is in the state~$s'$, we can write 
\begin{equation}
    \mathbb{P}(Y_{t} \leqslant M \,|\, B(0)\!=\!s) 
    = \sum_{s' \in S}  \mathbb{P}(B(t)\!=\!s'  \,|\, B(0)\!=\!s)\,\mathbb{P}(T[{s'}] \leqslant M\!-\!t)
    = {\boldsymbol p}_{t}^s\:\big(\,\mathbb{P}(T[{s'}] \leqslant M\!-\!t)\,\big)^\top_{s' \in S}.
    \label{condProbOffline}
\end{equation}
Hence, if we are able to determine the distribution of the travel time  $T[{s}]$ for each state $s \in S$, this would allow us to compute the conditional probability in \eqref{condProbOffline} for each $t$, which in turn would facilitate determining the maximizer~$t^\star $ in \eqref{optimal_departure_offline}. The evaluation of the travel-time distribution is outlined in detail in Section~\ref{sec:optdeptime}.

\medskip

\textit{Online Setting --- }
In the offline setting, the objective is to produce an optimal departure time $t^\star $ that depends solely on the current state of the background process $B(0)$. However, when a traveler waits for departure, new information
on the state of the background process becomes available. This in turn alters the distribution ${\boldsymbol p}_{t^\star }^s$ in~\eqref{distBOffline}. More specifically, for $u \leqslant t^\star $, the distribution of $B(t^\star )$ given $B(u) = s$ is 
\begin{equation}
\left(\,\mathbb{P}(B(t^\star ) = s' \,|\, B(u) = s)\,\right)_{s' \in S}
= \left(\,\mathbb{P}(B(t^\star -u) = s' \,|\, B(0) = s)\,\right)_{s' \in S}
= {\boldsymbol p}_{t^\star -u}^s \,e^{Q (t^\star -u)}.
\label{state_distribution}
\end{equation}
Since the distribution of $B(t^\star )$ changes as time $u$ progresses from $0$ to $t^\star  $, the on-time arrival probability of the driver changes as well. Ideally,  as a driver is waiting to depart,  their optimal departure time is updated such that it incorporates the latest state of the network.
%, and consequently, such that their on-time arrival probability remains unaltered. 
This way, as time passes, the traveler can request a new departure time that is given by
\begin{equation}
	t_u^\star  := \sup\big\{t \geqslant  u:\, \mathbb{P}\big(Y_{t} \leqslant M \mid B(u)\big) \geqslant  \eta\big\}.
\label{optimal_departure_online}
\end{equation}
Just as in the offline setting, $Y_{t}$ should be interpreted as the  arrival time, with the current time now being the request time $u$. Therefore, $t_u^\star $ is the latest departure time such that at time $u$ the on-time arrival probability is at least $\eta$, and this online departure time coincides with the offline departure time if $u = 0$ (conditional on $B(u)$ applying at time $0$). 

Again, just like in the offline setting, the on-time arrival probability at time of departure may differ from $\eta$. However, the request time will approach the departure time if a driver keeps updating their optimal departure time. Therefore, the on-time arrival probability at time of request will approach the on-time arrival probability at time of departure. Hence, in the online setting, it is in fact possible to find the latest departure time such that, on departure, a certain on-time probability can be satisfied.

\section{Deriving the Optimal Departure Time} \label{sec:optdeptime}
Recall that we consider the setting in which 
a vehicle wants to know its optimal
departure time for the traversal of an OD-pair in a network
$G = (N,A)$, in which vehicle speeds are described by the Markovian Velocity
Model that was presented in Section~\ref{subsec: Markov}. 
In Section~\ref{subsec:traveltimedist} we point out how the travel-time
distribution can be numerically evaluated by applying a discretization
procedure, while in Section~\ref{subsec:granularity} we discuss how the 
corresponding granularity can be determined.
{Then, Section~\ref{subsec:online} 
outlines how to compute the optimal departure time
for a given path, or, if only the destination
of the trip is known, how to obtain both
the optimal departure time and corresponding path to travel.
We consider the case 
in which the traveler only requests a departure time
once, i.e., the offline setting, as well as
how to extend this procedure to
the case in which the
traveler is allowed to receive departure-time advice 
updates, i.e., the online setting.
{
Note that, even though this section only treats
the compact version of the {\sc mvm}, the optimal departure
time procedures extend to the case the Markov
process $B(t)$ describes more detailed congestion phenomena,
of which two examples are provided in Appendix~\ref{app:recurrent}.}
}

% Then, Section~\ref{subsec:online} 
% outlines how to obtain the optimal departure time
% in case the traveler only requests a departure time
% once, i.e., the offline setting, as well as
% how to extend this procedure to
% the case in which the
% traveler is allowed to receive departure time advice 
% updates, i.e., the online setting.
% Section~\ref{subsec:offline} shows how to obtain the 
% optimal departure time in case the traveler only
% requests this departure time once, i.e., the offline setting.
% Section~\ref{subsec:online} demonstrates how to extend
% these results in the online setting, in which the
% traveler is allowed to receive departure time advice updates.
% In order to solve both these problems we will need the travel
% time distribution of the vehicle, which we
% derive in Section~\ref{subsec:traveltimedist}.

\begin{comment}
\begin{itemize}
    \item Via LS-transform, shortly explain why does not work
    \item Discretization
    \item Granularity: rule + figures to show resembles distribution well
\end{itemize}
\end{comment}

\subsection{Travel-Time Distribution} \label{subsec:traveltimedist}~\\
Consider a vehicle that departs at $t=0$ to traverse
a given path in the network $G = (N,A)$. The travel-time
distribution of this vehicle is completely determined by 
the velocity dynamics, described through the dynamics of the
Markovian background process~$B(t)$ and its initial state~$B(0)$. 
As both exact methods and Laplace inversion
fail (see Remark~\ref{remark:mixture} below), 
we use a discretization procedure to obtain an accurate approximation to
the travel-time distribution. Specifically, we discretize the moments
at which there is a potential transition in the driveable
speed (which corresponds to a transition of the background process).

\begin{remark} \label{remark:mixture}
{\em {
While traveling on a given link, 
the background process, and thus the
driveable vehicle speed, can in principle have
infinitely many transitions. Therefore, a
closed form distribution function
of the per-link travel time is unknown. 
% With the travel time of a path being the convolution of
% dependent individual link travel times, 
% a closed-form distribution
% function for $T_0$, the vehicle travel time, is unknown
% as well.
Since the Laplace-Stieltjes transform 
\textit{is} known \citep{Levering2021AConditions},
a natural procedure for obtaining the per-link travel-time distribution would be 
to rely on the numerical inversion of this transform.
Unfortunately, application of
common inversion methods (e.g., 
\cite{abate1992fourier}, \cite{iseger2009laplace}, or the
saddlepoint approximation \citep{Butler2007SaddlepointApplications}),
is not applicable. This is due to the fact that the per-link 
travel-time distribution is neither discrete nor continuous.
To see this, we consider
link $a$ with length denoted $d_{a}$, which takes
state $X_a(0) = s$ upon departure.
Now, there is a positive probability
that the state of $X_a(t)$ does not change while traveling link $a$,
and thus, that the link travel time equals $d_{a}/v_{a}(s)$.
However, as is easily seen, on the remainder of the domain
the per-link travel time has a continuous density.
}} $\hfill\Diamond$
\end{remark}

An easy fix for the challenges discussed in Remark~\ref{remark:mixture} would be to assume that the driveable vehicle speed on a link is fixed upon entering a link, in that is completely determined by the background state  upon entering (i.e.,  not affected by any transitions of the background process while  traversing the link). Doing so, the link travel time would reduce to a discrete distribution with travel-time probabilities that can easily be computed. While this procedure clearly gives a decent approximation for short links, it could perform poorly for  long links. This inspired us to the following idea: instead of assuming a fixed driveable vehicle speed for an entire link, we only assume fixed velocities for a certain (short) time interval. 
We now provide a detailed description of the procedure to obtain an approximation for the per-link travel-time distribution (which is also summarized in Algorithm~\ref{alg:distrsinglelink} below).

As a first step, we focus on the travel-time 
distribution for the traversal of a single link $a \in A$.
On this link, the driveable speed level
is either $v_1 := v_{a}(1)$ or $v_2 := v_{a}(2)$ for any $t \geqslant  0$
(see Section~\ref{subsec: Markov}).
Since we are solely given the state of $X_{a}(t)$ at $t = 0$, 
only the driveable speed upon departure is known.
Now, to be able to compute the travel-time distribution,
we discretize the moments at which the background process,
and consequently, the speed level, can change.

\begin{figure}[ht]
\centering
  \begin{subfigure}[b]{0.35\textwidth}
\begin{centering}
\begin{tikzpicture}[shorten >=1pt,node distance=2cm,on grid,auto]
    \node[state] (q_0) {$1$};
    \node[state] (q_1) [right=2 of q_0] {$2$};

    \path[->]
    (q_0) edge [bend left] node {$\alpha$} (q_1)
    (q_1) edge [bend left] node {$\beta$} (q_0);
\end{tikzpicture}
    \caption{\doublespacing Continuous.}
    \label{fig:discretizationXCont}
      \end{centering}
  \end{subfigure}
  \quad
  \begin{subfigure}[b]{0.35\textwidth}
  \begin{centering}
\begin{tikzpicture}[shorten >=1pt,node distance=2cm,on grid,auto]
    \node[state] (q_3) {$1$};
    \node[state] (q_4) [right=2 of q_3] {$2$};

    \path[->]
    (q_3) edge [loop left] node {$e^{-\alpha \delta}$} (q_3)
    (q_4) edge [loop right] node {$e^{-\beta \delta}$} (q_4)
    (q_3) edge [bend left] node {$1\!-\!e^{-\alpha \delta}$} (q_4)
    (q_4) edge [bend left] node {$1\!-\!e^{-\beta \delta}$} (q_3);
\end{tikzpicture}
    \caption{\doublespacing Discrete.}
    \label{fig:discretizationXDisc}
\end{centering}
  \end{subfigure}
    \caption{\doublespacing Transition rate diagrams for the continuous process $X_{a}(t)$ and its discretized version $X_{a}'(t)$.}
    \label{fig:discretizationX}
\end{figure}
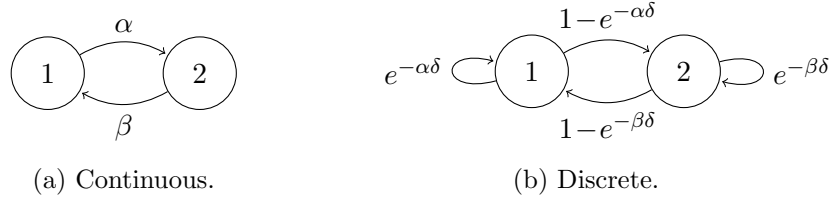

Concretely, given a (typically small) $\delta \in \mathbb{R}_{> 0}$, 
we define the Markov chain $X_{a}'(t)$
as a discrete-time version of $X_{a}(t)$,
at times $t=0,\delta,2\delta,\dots$.
That is, we set 
$X_{a}'(0) = X_{a}(0)$, and let the discrete-time Markov chain
$X_{a}'(t)$ (with $t=0,\delta,2\delta,\dots$) evolve
according to the diagram of Figure~\ref{fig:discretizationX}b.
Thus, if the process $X'_a(\cdot)$ is in state 1 at some time $m\delta$ (with $m\in{\mathbb N}_0$),
it is still there at $(m+1)\delta$ with probability $p := e^{-\alpha \delta}$, i.e., the probability that the continuous-time Markov chain $X_{a}(\cdot)$ does not jump to state 2 in a time interval of length $\delta$.
Alternatively, the process jumps to state 2 with
probability $1\!-\!p$,
i.e., the probability that
$X_{a}(\cdot)$ jumps
in a time interval of length $\delta$.
Note that, even though this latter 
probability incorporates the event
of multiple transitions of the continuous-time process $X_{a}(\cdot)$, 
our procedure is justified by the fact that,
if $\delta$ is chosen sufficiently small, 
the probability of more than one such
transition within a time interval
of length $\delta$ is $o(\delta)$ and therefore negligible; further
details regarding the choice of $\delta$ are discussed
in Section~\ref{subsec:granularity}.
In a similar fashion, if $X_{a}'(\cdot) = 2$,
the process stays in state~2 with probability $q := e^{-\beta \delta}$,
and moves to state 1 with probability $1\!-\!q$.
Now that we have described the dynamics of $X_{a}'(\cdot)$, 
the corresponding velocities can be defined: if, for $m \in \mathbb{N}_0$,
$X_{a}'(m\delta) = s$,
the speed level during the interval $[m\delta, (m\!+\!1)\delta)$
equals $v_{a}(s)$. Hence, during such an interval
of length $\delta$,
the speed level is constant.

Note that with the described velocity dynamics,
we are able to iteratively compute the travel time
distribution on link~$a$.
That is, consider the case that upon departure,
link~$a$ is free of incidents, i.e., $X_{a}'(0) = 1$
(a similar procedure can be followed for
the alternative situation).
Then, with a constant speed~$v_1$ during $[0,\delta)$, 
the traveled distance at time~$\delta$ equals 
$v_1\delta$ with probability~1.
Now, the traveled distance 
at time $2\delta$ equals
$(v_1\!+\!v_2)\delta$ with probability $1\!-\!p$ 
(the probability that $X_{a}'(\delta)\!=\!2$),
and $2v_1\delta$ with probability $p$ (the probability
that $X_{a}'(\delta)\!=\!1$).
These two cases each lead to two potential travel times
at time $3\delta$:
in case $X_{a}'(\delta) = 1$, the traveled
    distance equals $3v_1\delta$ with probability $p^2$,
    and equals $(2v_1\!+\!v_2)\delta$ with probability $p(1\!-\!p)$, and in case $X_{a}'(\delta) = 2$, the traveled
    distance equals $(2v_1\!+\!v_2)\delta$ with probability $(1\!-\!p)(1\!-\!q)$,
    and $(v_1\!+\!2v_2)\delta$ with probability $(1\!-\!p)q$.
We can iteratively continue these computations,
in which every (state, distance, probability)-tuple at $t = m\delta$
generates two tuples for $t = (m\!+\!1)\delta$. Thus, any
tuple in $t = (m\!+\!1)\delta$ has a so-called \textit{ancestor}
in $t = m\delta$.

To obtain the travel-time
distribution of this link, recall that $d_{a}$ is the distance of the link,
and thus the total distance to travel.
Therefore, if at $t = m\delta$ a tuple ($s_0, d_0, p_0$) with 
traveled distance $d_0 < d_{a}$
generates a tuple ($s_1, d_1, p_1$) for which $d_1 \geqslant  d_{a}$, 
then there is a probability $p_1$
that the travel time of the link equals 
$m\delta + (d_{a}\!-\!d_0)/v_{a}(s_1)$.
Namely, in this scenario, 
after $m\delta$ time the vehicle
still needs to traverse a distance~\mbox{$d_{a}\!-\!d_0$}, which it travels
with the constant speed level belonging to state $s_1$.
The collection of travel-time values and 
corresponding probabilities
that are iteratively found in this fashion, 
form the travel-time distribution. We  denote this collection for
link $a_i$ given $X_{a_i}(0) = j$ as $L_{a_i}^j$.
Observe that the tuples ($s, d, p$) for which
$d \geqslant  d_{a}$ do not serve as an ancestor in a new iteration, such that
the new iteration only continues with tuples
for which the traveled distance has not yet exceeded the 
length of the link.
Moreover, since the traveled distance grows
with $v_1\delta$ or $v_2\delta$ every step, there are only
finitely many steps before $d_{a}$ is exceeded,
and consequently, the iterative procedure terminates in a finite number of steps.

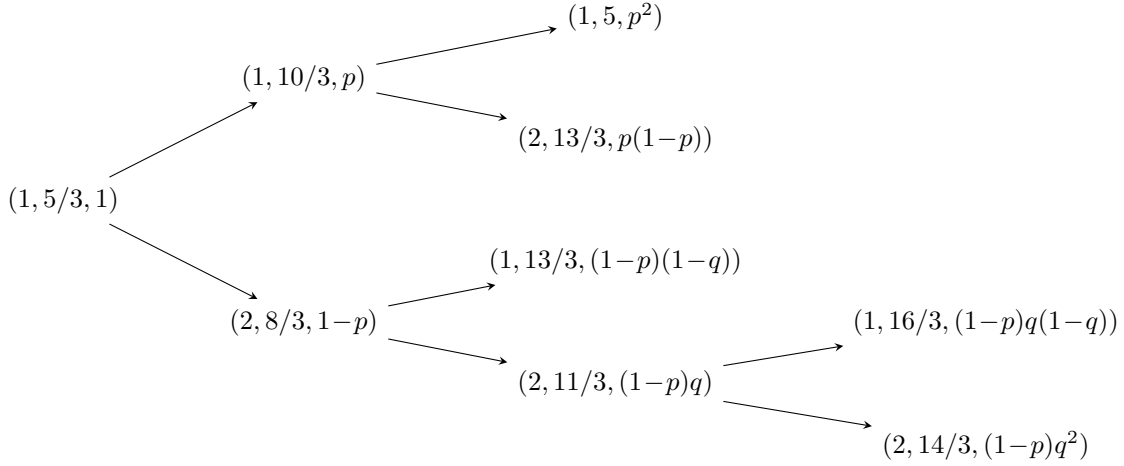
\begin{figure}[ht]
\begin{centering}
\begin{tikzpicture}[>=stealth,sloped]
    \matrix (tree) [%
      matrix of nodes,
      minimum size=0.5cm,
      column sep=1.2cm,
      row sep=0.2cm,
    ]
    {
          &   & \small{($1,5,p^2$)} &  \\
          & \small{($1,10/3,p$)} &   &  & \\
          &   & \small{($2,13/3,p(1\!-\!p)$)} & \\
        \small{($1,5/3,1$)} &   &   &  & \\
          &   & \small{($1,13/3,(1\!-\!p)(1\!-\!q)$)} & \\
          & \small{($2,8/3,1\!-\!p$)} &   & \small{($1,16/3,(1\!-\!p)q(1\!-\!q)$)} \\
          &   & \small{($2,11/3,(1\!-\!p)q$)} &  & \\
          &   &   & \small{($2,14/3,(1\!-\!p)q^2$)} \\
    };
    \draw[->] (tree-4-1) -- (tree-2-2) node [midway,above] {};
    \draw[->] (tree-4-1) -- (tree-6-2) node [midway,below] {};
    \draw[->] (tree-2-2) -- (tree-1-3) node [midway,above] {};
    \draw[->] (tree-2-2) -- (tree-3-3) node [midway,below] {};
    \draw[->] (tree-6-2) -- (tree-5-3) node [midway,above] {};
    \draw[->] (tree-6-2) -- (tree-7-3) node [midway,below] {};
    \draw[->] (tree-7-3) -- (tree-6-4) node [midway,above] {};
    \draw[->] (tree-7-3) -- (tree-8-4) node [midway,below] {};
  \end{tikzpicture}
  \end{centering}
\caption{\doublespacing Tree-structure iterative procedure from 
Example~\ref{ex:tree}.}
\label{fig:tree}
  \end{figure}

% \begin{figure}[ht]
% \begin{centering}
% \begin{forest}
% for tree={
%       l sep = 30pt,
%       if n children=0{
%       font=\itshape,
%       tier=terminal,
%     }{},
%   }
% % for tree = {
% % l sep = 30pt,
% % }
%   [{\small{($1,5/3,1$)}}
%   [{\small{($1,10/3,p$)}}
%      [{\small{($1,5,p^2$)}}
%      [{\small{\textit{1/25h w.p. $p^2$}}}]
%      ]
%      [{\small{($2,13/3,pp'$)}}
%      [{\small{\textit{2/45h w.p. $pp'$}}}]
%      ]
%     ]
%     [{\small{($2,8/3,p'$)}}
%       [{\small{($1,13/3,p'q'$)}}
%       [{\small{\textit{7/150h w.p. $p'q'$}}}]
%       ]
%       [{\small{($2,11/3,p'q$)}}
%         [{\small{($1,16/3,p'qq'$)}}
%         [{\small{\textit{4/75h w.p. $p'qq'$}}}]
%         ]
%         [{\small{($2,14/3,p'q^2$)}}
%         [{\small{\textit{1/18h w.p. $p'q^2$}}}]
%         ]
%       ]
%     ]
%   ]
% \end{forest}
% \end{centering}
% \caption{Tree-structure iterative procedure from 
% Example~\ref{ex:tree}}
% \label{fig:tree}
% \end{figure}

% \begin{tikzpicture}
 
% \node {$1,5/3,1$} [sibling distance = 2.5cm]
%     child {node {child 1} 
%     child {node {grandchild 1}}
%     child {node {grandchild 2}}} 
%     child {node {child 2}
%     child {node {grandchild 1}}
%     child {node {grandchild 2}}};
 
% \end{tikzpicture}

\begin{algorithm}[ht]
 \KwResult{\doublespacing list $L$ of (travel time, probability) values.}
 \doublespacing
 Notation: for transparency, omitted subscripts for $v_{a}(s)$,
 $X_{a}(s)$ and $d_{a}$\;
 Given: $\delta \in \mathbb{R}_{>0}, X(0) = s, p = e^{-\alpha \delta}, q = e^{-\beta \delta}$\;
 Initialization:  $L = \emptyset$, $S = \{(s,\delta v(s),1)\}$, $S' = \emptyset$, $m = 0$\;
 \While{$S$ nonempty}{
 1. \ForEach{$(s',d',p')$ in $S$}{
 a. Extract $(s',d',p')$ from $S$\;
 b. Compute $p_1 = p'(1\{s'\!=\!1\}p + 1\{s'\!=\!2\}(1\!-\!q))$ and $p_2 = p'(1\{s'\!=\!1\}(1\!-\!p) + 1\{s'\!=\!2\}q)$\;
 c. \textbf{If} $d'+v(1)\delta < d$ \textbf{then} 
    append $(1,d'+v(1)\delta,p_1)$ to $S'$. \textbf{Else}
    append $(m\delta + (d\!-\!d')/v(1),p_1)$ to $L$\;
 d. \textbf{If} $d'+v(2)\delta < d$ \textbf{then} 
    append $(2,d'+v(2)\delta,p_2)$ to $S'$. \textbf{Else}
    append $(m\delta + (d\!-\!d')/v(2),p_2)$ to $L$\;    
}
 2. $m \mathrel{+}= 1$ \;
 3. Set $S = S', S' = \emptyset$;
 }
 \caption{Travel-time distribution of a single link $a$}
 \label{alg:distrsinglelink}
\end{algorithm}

\begin{example} \label{ex:tree}
{\em {
Figure~\ref{fig:tree} displays the iterative steps of the
procedure via a tree structure, for obtaining
the travel-time distribution on a link $a \in A$
with $d_a = 4$ km, $v_a(1) = 100$ km/h, $v_a(2) = 60$ km/h,
and $X_{a}(0) = 1$. In the discretization,
we only allow speed transitions at full minutes, i.e.,
we let $\delta = 1/60$ h.
Iteration at a branch of the tree stops in case
the traveled distance exceeds 4 km, and results
in a travel-time value.
For example, the upper branch yields a travel time
of $2/60 + (4\!-\!10/3)/100 = 1/25$ h, and the branch
directly below a
travel time of $2/60 + (4\!-\!10/3)/60 = 2/45$ h.
}}\hfill$\Diamond$
\end{example}

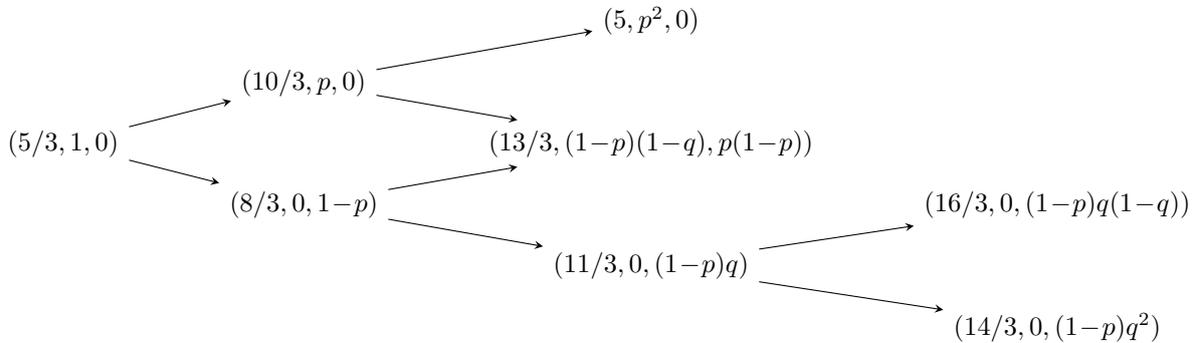
\begin{figure}[ht]
\begin{centering}
\begin{tikzpicture}[>=stealth,sloped]
    \matrix (tree) [%
      matrix of nodes,
      minimum size=0.5cm,
      column sep=1.2cm,
      row sep=0.2cm,
    ]
    {
          &   & \small{($5,p^2,0$)} &   \\
          & \small{($10/3,p,0$)} &   &   \\
        \small{($5/3,1,0$)} &   &  \small{($13/3,(1\!-\!p)(1\!-\!q),p(1\!-\!p)$)} &   \\
          & \small{($8/3,0,1\!-\!p$)} &   & \small{($16/3,0,(1\!-\!p)q(1\!-\!q)$)} \\
          &   & \small{($11/3,0,(1\!-\!p)q$)} &   \\
          &   &   & \small{($14/3,0,(1\!-\!p)q^2$)} \\
    };
    \draw[->] (tree-3-1) -- (tree-2-2) node [midway,above] {};
    \draw[->] (tree-3-1) -- (tree-4-2) node [midway,below] {};
    \draw[->] (tree-2-2) -- (tree-1-3) node [midway,above] {};
    \draw[->] (tree-2-2) -- (tree-3-3) node [midway,below] {};
    \draw[->] (tree-4-2) -- (tree-5-3) node [midway,above] {};
    \draw[->] (tree-4-2) -- (tree-3-3) node [midway,below] {};
    \draw[->] (tree-5-3) -- (tree-6-4) node [midway,above] {};
    \draw[->] (tree-5-3) -- (tree-4-4) node [midway,below] {};
  \end{tikzpicture}
  \end{centering}
  \caption{\doublespacing Binomial tree representation corresponding to 
Example~\ref{ex:tree}.}
  \label{fig:bintree}
  \end{figure}

\begin{remark}
{\em {
We observe that, since we are only working with
two constant speed values per time step, 
the iterative procedure can, in this
special case, also be represented by a binomial tree. 
To be able to generate such a tree and to compute the
resulting travel-time distribution, we look,
contrary to Algorithm~\ref{alg:distrsinglelink},
at tuples of the form ($d, p_1, p_2$), with $d$ is the traveled distance,
and $p_1$ ($p_2$, respectively) the probability of the scenario in which
the vehicle had speed $v_1$ ($v_2$, respectively) in the last time step.
% and $p_2$ the probability of the scenario in which
% the vehicle had speed $v_2$ in the last time step .
Note that we need to separate these two probabilities,
as the corresponding two scenarios affect the probabilities of the next
time step differently. 
Figure~\ref{fig:bintree} shows the binomial
tree corresponding to Figure~\ref{fig:tree}.

Importantly, a binomial tree grows by
maximally one item per time step,
making this procedure particularly efficient.
Indeed, after $m$ time steps, the 
number of ancestors in a binomial
tree is only of 
the order $m$, whereas  
the number of ancestors in Algorithm~\ref{alg:distrsinglelink}
would be of order $2^m$.
}} \hfill$\Diamond$
\end{remark}

% \begin{remark}
% % Only binomial tree if indeed incidents
% %at boundary of interval. Middle of expectation would
% %  yield much larger trees.
% {\em {
% In the example above, the iteration tree is very small,
% with just four layers. 
% If the number of required layers is much higher, the iteration tree
% becomes much larger, as, typically, every tuple 
% generates two new tuples. Even though this
% makes computational complexity of the iteration procedure very high,
% the iterative procedure \textit{can} be performed 
% efficiently. That is, since every time step a traveled
% distance is either increased by $v_1x$ or $v_2x$,
% purely focusing on the traveled distances would yield
% a binomial tree. Importantly, such a tree grows
% maximally one item per time step.

% To be able to generate such a tree and compute the
% resulting travel time distribution, we look
% at tuples of the form ($d, p_1, p_2$),
% in which $d$ is the traveled distance,
% $p_1$ the probability of the scenario in which
% the last time step the vehicle had a speed $v_1$,
% $p_2$ the probability of the scenario in which
% the last time step the vehicle had a speed $v_2$.
% We need to separate these two probabilities,
% as they affect the probabilities of the next
% time step differently. 
% Figure~\ref{fig:bintree} shows the binomial
% tree corresponding to Figure~\ref{fig:tree}.
% }} \hfill$\Diamond$
% \end{remark}

Now, knowing how to compute the travel-time distribution for the traversal
of an individual link, we may redirect our focus to the travel-time distribution of a full path $P$.
Without loss of generality, we let
$P = \{a_1,\dots,a_m\}$.
Writing $s_i$ for the initial state 
$X_{a_i}(0)$, the travel-time distribution
for $a_{1}$ is known and given by the list $L_{a_1}^{s_1}$,
consisting of pairs ($t,p$), with $t$ a travel-time value and $p$ the corresponding probability.
Now, let us have a list~$L$ of $(t,p)$ pairs that form
the travel-time distribution on the subpath $\{a_1,\dots,a_i\}$.
Then, we observe that, if the travel time for $\{a_1,\dots,a_{i}\}$
equals $t_1$, the probability that $X_{a_{i+1}} = j$
upon entering $a_{i+1}$ is the ($s_{i+1}, j$)-th
index of $e^{t_1Q_{a_{i+1}}}$, with $Q_{a_{i+1}}$ 
the transition rate
matrix of $a_{i+1}$ as defined in \eqref{eq:qmatrix}.
Therefore, the travel-time distribution~$T[s]$ for traversing
$\{a_1,\dots,a_{i+1}\}$ with initial state~$s$ is given by the list 
\begin{equation}
    \big\{\big(t_1 + t_2, p_1 \cdot p_2 \cdot [e^{t_1Q_{a_{i+1}}}]_{(s_{i+1}, j)}
    \big) \,|\, (t_1,p_1) \in L, 
    (t_2,p_2) \in L_{a_{i+1}}^j, j=1,2 \big\}.
    \label{eq:pathdist}
\end{equation}
Notably, by iteratively setting $i = 1,\dots,m\!-\!1$,
we obtain the travel-time distribution of the path $P$. 
To prevent that the number of elements in the list $L$ 
grows exponentially with the number of links in the path,
we aggregate the travel-time values into equally sized 
bins after every iteration.

\subsection{Granularity} \label{subsec:granularity}~\\
Now that we have outlined the discretization procedure that allows us to approximate the per-path travel-time distribution, it remains to choose an appropriate value for the length of the time intervals $\delta$. The choice of $\delta$ affects the number of steps in the iterative procedure, which we call the granularity. It is clear that choosing a relatively large value of $\delta$ renders the algorithm fast. However, for two reasons the value of $\delta$ may not be chosen too large, as we point out now.

First, for a given link $a \in A$, we used in Section~\ref{subsec:traveltimedist} the approximation 
$$\mathbb{P}(X_{a}(\delta) = 1 \mid X_{a}(0) = 1) \approx p := e^{-\alpha \delta}.$$ Note that $e^{-\alpha \delta}$ is the probability that $X_{a}(0)$ has no transitions in the time interval $\delta$. However, the probability $\mathbb{P}(X_{a}(\delta) = 1 \mid X_{a}(0) = 1)$ also contains the events in which $X_{a}(\cdot)$ jumps multiple times, and returns to its original state to remain in that state until the end of the time interval $\delta$. Therefore, in order to ensure that the approximations we used in Section~\ref{subsec:traveltimedist} are justifiable, we require the probability of two or more transitions in the time interval $\delta$ to be negligible.
With $E(\gamma)$ an exponentially distributed random variable with parameter $\gamma$ and $f_\gamma(\cdot)$ its density, the probability of two or more transitions conditional on $X_a(0)=1$ is given by
\begin{align*}\int_0^\delta f_\alpha(t) \int_0^{t-\delta} f_\beta(s)\,{\rm d}s\,{\rm d}t &= \int_0^\delta \alpha e^{-\alpha t}(1-e^{-\beta(\delta-t)})\,{\rm d}t\\&=1-e^{-\alpha\delta}-\frac{\alpha}{\alpha-\beta}\big(e^{-\beta\delta}-e^{-\alpha\delta}\big)=1-\frac{\alpha e^{-\beta\delta}-\beta e^{-\alpha\delta}}{\alpha-\beta}.\end{align*}
Due to the fact that this expression is symmetric in $\alpha$ and $\beta$, 
it also equals the probability of two or more transitions conditional on $X_a(0)=2$. 
This concretely means that, for a given small value of $\varepsilon_1 > 0$ (for instance~$0.01$), the interval length $\delta$ should be chosen sufficiently small so that
\[
%\max\Big\{\mathbb{P}(), \mathbb{P}()\Big\}
1-\frac{\alpha e^{-\beta\delta}-\beta e^{-\alpha\delta}}{\alpha-\beta} < \varepsilon_1.
\addtag
\label{eq:granerror1}
\] 
%Note that we have a maximum since we need to consider both the case $X_{a}(0) = 1$ and $X_{a}(0) = 2$. 
As $\delta$ should be chosen small, we can {simply} evaluate \eqref{eq:granerror1} working with its second-order Taylor approximations in $\delta$ at $0$. Indeed, as $\delta\downarrow 0$,
\[1-\frac{\alpha e^{-\beta\delta}-\beta e^{-\alpha\delta}}{\alpha-\beta}  = \frac{\frac{1}{2}\alpha^2\beta\delta^2-\frac{1}{2}\alpha\beta^2\delta^2}{\alpha-\beta}+o(\delta^2)=\frac{1}{2}\alpha\beta\delta^2+o(\delta^2).\]
We thus find that \eqref{eq:granerror1} reduces to
\[
\delta < \left(\frac{2\varepsilon_1}{\alpha\beta}\right)^{1/2}.
\addtag
\label{eq:granerrorapprox1}
\]
From \eqref{eq:granerrorapprox1}, it is clear that $\delta$ should be chosen smaller for higher transition rates $\alpha, \beta$. As higher transition rates result in a higher probability of the occurrence of two or more transitions within a fixed time interval, it is also intuitively clear that $\delta$ should be decreasing in $\alpha, \beta$ in order to constrain this approximation error. 

Second, we used that the driveable vehicle speed on a link is fixed for the time period $\delta$. If, however, a transition from, say, state 1 to state 2 occurs at time $t<\delta$, the actual traversed distance in the time interval $\delta$ would be $v_1t + v_2(\delta-t)$ as opposed to
our approximated $v_1\delta$. This approximation error is more substantial if the difference between the velocities $v_1$ and $v_2$ increases. In case $X_{a}(0) = 1$, we thus want to choose $\delta$ sufficiently small such that for a chosen small $\varepsilon_2 >0$,  the expected error is bounded:
    \begin{equation*}
        (1-e^{-\alpha \delta})\left|\,\int_0^\delta 
        \frac{\alpha e^{-\alpha t}}{1-e^{-\alpha \delta}}
        (v_1t + v_2(\delta-t))\,{\rm d}t - v_1\delta\,
        \right|< d \cdot \varepsilon_2.
    \end{equation*}
This can be rewritten to
    \begin{equation}
        \big|v_1-v_2\big|\,\left(\frac{1}{\alpha}(1-e^{-\alpha \delta})-\delta\right)  < d \cdot \varepsilon_2.
        \label{eq:granerror2}
    \end{equation}
A similar condition can be obtained for the case $X_{a}(0) = 2$, but with the $\alpha$ replaced by $\beta$ in \eqref{eq:granerror2}. Note that the error $\varepsilon_2$ is multiplied by the length $d$, since we are interested in the error relative to the length of a link. Just as we did for the first approximation error, we can simplify \eqref{eq:granerror2} further by considering its Taylor approximations in $\delta$ at $0$. Doing this for both  $X_{a}(0) = 1$ and $X_{a}(0) = 2$, we find
\[
\delta < \left(\frac{2\varepsilon_2}{\frac{|v_1-v_2|}{d}\alpha}\right)^{1/2} \quad \text{and} \quad \delta < \left(\frac{2\varepsilon_2}{\frac{|v_1-v_2|}{d}\beta}\right)^{1/2}.
\addtag
\label{eq:granerrorapprox2}
\]
Combining conditions \eqref{eq:granerrorapprox1} and \eqref{eq:granerrorapprox2}, we conclude that both the probability of two or more transitions and the expected error are sufficiently small, whenever we pick
\[
\delta < \min\left\{\left(\frac{2\varepsilon_1}{\alpha\beta}\right)^{1/2},\left(\frac{2\varepsilon_2}{\frac{|v_1-v_2|}{d}\max\{\alpha,\beta\}}\right)^{1/2}\right\}.
\addtag
\label{eq:granerrorapproxcombined}
\]
\subsection{Optimal Departure Time} \label{subsec:online}~\\
%With the discretization procedure outlined in Sections~\ref{subsec:traveltimedist}
%and~\ref{subsec:granularity} 
So far, we have described how to compute the travel-time
distribution $T[s]$ for a vehicle traversing a path 
$P$ and departing when $B(t)$ is in state $s$,
for any $s \in S$.
Using knowledge of the dynamics of the background process~$B(t)$, this allows us to 
compute the on-time arrival probability for 
a path $P$ and any departure time~\mbox{$t \geqslant  0$}.
Observing that the on-time arrival probability is monotonically decreasing
in the departure time, by performing an elementary bisection we can determine the optimal (the latest, that is)
time to depart on path $P$ for a given on-time arrival probability.
This monotonicity can moreover be used to extend the
results to the case in which the user only specifies
the origin and destination, instead of the specific path to travel.
The procedure then compares the departure times
for different paths, selects the latest, and outputs both $t^\star $ and
the corresponding path to travel.

% For departure at $t = 0$ we are able to compute the travel time distribution, or equivalently,
% the arrival time distribution, of
% a vehicle traversing a path in the network. Now, using the dynamics of $B(t)$ in $[0,t']$, 
% we can use the same techniques to compute the arrival time distribution for 
% any departure time $t' > 0$.
% Thus, for any departure time $t' \geqslant  0$ we are able to compute the on-time
% arrival probability. By observing that the on-time arrival probability is monotone decreasing
% in the departure time, use of the bisection algorithm will yield the optimal
% departure time for a given on-time arrival probability.

\subsubsection{Optimal departure time for a path --- offline setting}~\label{subsec:pathoffline}
One natural way to obtain the on-time arrival probability on
a path $P$ for a departure
time $t \geqslant  0$ was presented
in Section~\ref{subsec: departureproblem}, namely, computing the product in \eqref{condProbOffline}.
To this end, we need to evaluate $\mathbb{P}(T[{s'}] \leqslant M\!-\!t)$ for all $s' \in S$.
Note that the distribution of $T[{s'}]$ can be derived
by the presented discretization procedure, which outputs a list of 
(travel time, probability)-pairs. 
By summing all probabilities for which the corresponding
travel-time value does not exceed $M\!-\!t$, we obtain a value for 
$\mathbb{P}(T[{s'}] \leqslant M\!-\!t)$.

However, computing the arrival probability
via \eqref{condProbOffline} is typically time consuming, as 
it requires us to compute the distribution of $T[s]$ for all $s \in S$.
Fortunately, there is an alternative procedure for which only one travel-time distribution needs to be derived.
That is, realize that departing at time $t \geqslant  0$ to
traverse path $P$ can be viewed as departing at time $0$ and,
before entering $P$, first a fictional link has to be traversed
for which the travel time
equals $t$ with probability~1.
Therefore, we directly obtain the 
distribution of the arrival time $Y_{t}$ by computing the travel
time distribution of this partly fictional path
(with departure at time 0). Then,
the on-time arrival probability can be found
by summing all probabilities for which the
corresponding travel time does not exceed~$M$.

Now, since we are able to compute the on-time arrival
probability on the path $P$ for a given departure time $t$ and, moreover,
this probability is monotonically decreasing in the departure time,
we can use bisection to find the optimal departure time 
for the given on-time arrival probability $\eta$. The monotonicity
of the on-time probability follows directly from
Proposition~1 of \cite{Levering2021AConditions}, which gives
that for a path consisting of a single link, $t' \leqslant t$ implies
\begin{equation*}
    t' + T_{t'} \leqslant t + T_{t},
\end{equation*}
and therefore, $Y_{t'} \leqslant Y_t$.
Since the minimum and maximum velocity on all links of the path
are known, the minimum travel time $t_{\text{min}}$ and 
the maximum travel time $t_{\text{max}}$ are known as well,
so that
\begin{equation*}
    I_0 := [\max\{0,M-t_{\text{max}}\}, \max\{0,M-t_{\text{min}}\}]
\end{equation*} 
serves as natural starting 
interval for the bisection method.
First, we check if $I_0$ equals $[0,0]$,
since in that case, the minimum travel time is at least $M$,
and $t^\star  := -\infty$.
Second, we check the on-time arrival probability at the
left boundary. If this probability is below $\eta$,
$t^\star  := -\infty$ as well.
%the departure advice is to leave at $t = 0$ as well.
In case neither is true, we apply the bisection algorithm until
we obtain the latest departure time for which the on-time
probability is at least $\eta$ (which is guaranteed to exist
by the first two checks).

Algorithm~\ref{alg:offline} now summarizes the complete procedure
for obtaining the optimal departure time~$t^\star $ to traverse  
a path $P$ for a given on-time arrival 
probability $\eta$ in the offline setting,
in which a traveler requests the value of $t^\star $ once.
It is important to note that we can precompute the travel-time distribution for every link $a_i$ in the network, and
for every initial state $s$ of this link. We can then store these distributions
as the lists $L_{a_i}^s$. Then, upon an optimal departure time
request of a vehicle, we can directly use these distributions,
and do not need to compute them on-the-spot.
Thus, the computational costs are only determined by
Part II of the algorithm.

\begin{algorithm}[ht]
 \KwResult{\doublespacing optimal departure time for traversing path $(a_1,\dots,a_m)$ within time $M$ with probability $\eta$, given $B(0) = s$.}
 \textcolor{white}{a}\\
 \doublespacing
 Part I (precomputations):\\
 \For{$i = 1,\dots,m$}{
    1. for $a_i$, compute $\delta$ via \eqref{eq:granerrorapproxcombined}\;
    2. Compute $L_{a_i}^1$ and $L_{a_i}^2$ with Algorithm~\ref{alg:distrsinglelink}
    and $\delta \in \mathbb{R}$ from step 1\;
 }
 \textcolor{white}{a}\\
 Part II (on-the-spot computations):\\
 1. Define the function:\\
 \SetKwProg{OnTimeProbability}{OnTimeProbability}{}{}
 \OnTimeProbability{$(t, s, M)$}{
 Set $L = \{(t,1)\}$\;
 \For{$i = 1,\dots,m$}{
 Set $L = \big\{\big(t_1 + t_2, p_1 \cdot p_2 
 \cdot [e^{Q_{a_{i}}t_1}]_{(s_{i}, j)}
 \big) \,|\, (t_1,p_1) \in L, 
 (t_2,p_2) \in L_{a_i}^{j}, j=1,2 \big\}$\;
 }
 Return $\sum_{(t_1,p_1) \in L : t_1 \leqslant M} p_1$\; 
 }
 2. Compute $t_{\text{min}} = \sum_{i=1}^m d_{a_i}/v_{a_i}(2)$
 and $t_{\text{max}} = \sum_{i=1}^m d_{a_i}/v_{a_i}(1)$\;
 3. Set $I_0 = [\max\{0,M-t_{\text{max}}\}, \max\{0,M-t_{\text{min}}\}]$\;
 4. \uIf{$I_0 = [0,0]$ or \em{\textbf{\text{OnTimeProbability}}}$(\max\{0,M-t_{\text{max}}\},s,M)    \leqslant \eta$}{Return $t^\star  = -\infty$\;}
    \uElse{Use the bisection method with initial interval $I_0$ on the function
    \textbf{OnTimeProbability}($t,s,M$), until obtain latest departure time for which
    output is at least $\eta$. Return this departure time $t^\star $.}
 \caption{Optimal departure time offline setting}
 \label{alg:offline}
\end{algorithm}

\subsubsection{Optimal departure time for an OD-pair --- offline setting} \label{subsub:od}
We can now extend the results to the natural setting
in which a user does not specify the complete path to travel, but only
the destination point $k^\star$. Then, the traveler does not just request
the optimal departure time to reach this endpoint with
on-time arrival probability~$\eta$, but also requests the specific path that guarantees
this probability $\eta$. 
% Note that we assume that the traveler
% follows this a-priori selected path, and 
% is not allowed to adapt the route while traveling. 
% We briefly discuss the more complex setting in which a 
% traveler is allowed to deviate from the path in Section~\ref{sec:concl}.
Note that, with $\mathcal{P}$ the set of all paths to the endpoint,
we can simply compute the optimal departure time for all
$P \in \mathcal{P}$, and output the path with the latest
departure time. However, since the size  $|\mathcal{P}|$ of such paths
is typically huge, the above procedure is 
not applicable in practical settings. 
Therefore, we present two alternative methods:
an exact procedure that is still somewhat computationally demanding,
and a very efficient, near-optimal method:

% The first method is based on a label-correcting procedure,
% and is guaranteed to output the optimal path and departure time,
% but at a computational complexity cost. 

\begin{itemize}
\item[$\circ$] Bisection method: similar to Algorithm~\ref{alg:offline}, the first procedure
uses a bisection algorithm to obtain the optimal
departure time and corresponding path. 
{That is, for a given departure
time~$t$, it uses a label-correcting algorithm
to output the path with maximum on-time arrival probability.
This algorithm (Algorithm~\ref{alg:od1}) will be
described in more detail below.
Now, as the maximum of monotonically decreasing functions, 
the on-time arrival probability outputted
by the label-correcting algorithm is again a monotonically
decreasing function of the departure time $t$.}
Consequently, bisection can indeed
be employed to find the optimal departure time
and corresponding path. 
Note that even though the bisection 
method is guaranteed to find the optimal
path and departure time, the computational complexity
of the label-correcting algorithm may prohibit practical
application.
\item[$\circ$] $k$-shortest path method: 
% even when only considering non-dominant paths,
% the label-correcting algorithm used in the method above
% may become prohibitively slow
% in large networks.
% Therefore, we propose an efficient alternative 
% for the computation of the optimal departure
% time for an OD-pair. In this
% second procedure, 
we simply compute 
the optimal departure
time for a small subset of $\mathcal{P}$.
Concretely, we compute the optimal departure time
for the subset of $k$ shortest paths 
(e.g.\ in distance) to the destination.
Thus, with $\mathcal{P}'$ the set of
$k$ shortest path as found via Yen's algorithm 
\citep{Yen1970, Yen1971}, we use Algorithm~\ref{alg:offline}
to compute the optimal departure time for all $p \in \mathcal{P}'$,
and output the path with latest departure time.
Note that, for $k$ large enough,
$\mathcal{P}'$ will typically contain the path which
yields the latest departure time. 
Importantly, as the optimal departure times for the
different paths can be computed in parallel, 
this method can indeed be employed in a highly efficient way.
\end{itemize}

We are left with describing the label-correcting algorithm, used within the bisection method,
to output the path 
with maximum on-time arrival probability
for a given departure time $t$.
This algorithm, outlined in Algorithm~\ref{alg:od1},
is an A$^\star$-algorithm 
in the same spirit as the algorithm 
presented in \cite{chen2014reliable},
assigning a label set to every node in
the graph and updating these sets iteratively.
The label set of a node is used to store
travel-time distributions of 
paths from the origin to that node,
when departing from the origin at time $t$.
Initially, the label set of the
origin consists of a single element, namely
the distribution $L_0 = \{(t,1)\}$, whereas
the other label sets start empty.
% Every iteration, a travel time distribution
% from one of the label sets is chosen,
% and used to compute travel time distributions
% to neighbors of the node corresponding to the label set.

The iteration uses a queue $q$, whose elements
are tuples of length four. Every tuple consists
of a node in the graph~($k$), a path from the
origin to this node~($P$), the list of (travel time, probability)-pairs
forming the travel-time distribution of this path~($L$),
and an upper bound of the maximum on-time arrival probability
for any path from origin to destination that has
$P$ as subpath~($\alpha$). 
This upper bound $\alpha$ is computed as 
$\mathbb{P}(T_{L} + t_{\text{min}}(k,k^\star) \leqslant M)$, in 
which $t_{\text{min}}(k,k^\star)$ is the minimum
travel time from $k$ to the destination $k^\star$
and $T_L$ a random variable with distribution $L$.
At the start, $q = \{\{\alpha,L_0,k_0,\text{path: }\{k_0\}\}\}$,
with $\alpha = 0$ if $t + t_{\text{min}}(k_0,k^\star) > M$
and $\alpha = 1$ otherwise. In case $\alpha = 0$,
all paths from $k_0$ to $k^\star$ 
will have a zero probability of on-time arrival,
thus the algorithm stops and outputs probability 0.
Otherwise, we continue. 
Now, every iteration step,
the element with minimum $\alpha$-value is extracted
from the queue. 
Then, for every neighbor $k'$ of
the node $k$ that is not in $P$, 
the travel-time distribution~$L'$ for traversing subsequently $P$ and
the link $(k,k')$ can be computed via \eqref{eq:pathdist}.
Clearly, neighbors $k' \in P$ are omitted as extending
$P$ with $(k,k')$ would yield a suboptimal path 
containing a loop.
After computing $\alpha'$, the upper bound
to the maximum on-time arrival probability for
the travel-time distribution $L'$,
we insert $L'$ into the label set of $k'$ and
the tuple $\{\alpha', L', k', P + k'\}$
into $q$. This ends the current iteration step,
and a new minimum element is extracted from $q$.
The algorithm is terminated
if the third element from the extracted tuple
from $q$ equals the destination $k^\star$. 

Now, to improve the speed of the procedure,
we perform one extra step before inserting a
new tuple in the label set of $k'$ and the queue $q$. 
That is, we check if the distribution $L'$ dominates
one of the distributions already stored in the label
set of $k'$, or
vice versa.
%alternatively,
%if there are distributions in the label
%set that dominate $L'$. 
In this setting, dominance refers to 
(first-order) stochastic dominance:
for two paths $p_1,p_2$ from $k_0$ to
$k'$ with respectively travel-time distributions $L_1$ and $L_2$, we say
that $L_1$ dominates $L_2$ if  
$\mathbb{P}(T_{L_2} \leqslant t) \leqslant \mathbb{P}(T_{L_1} \leqslant t)$ 
for all $t > 0$ and $\mathbb{P}(T_{L_2} \leqslant t) < \mathbb{P}(T_{L_1} \leqslant t)$ for at least one $t > 0$.
Now, if $L_1$ would indeed dominate a distribution
$L_2$ in the label set of $k'$,
then, for every path $p'$ from $k'$ to
the destination~$k^\star$, 
the on-time arrival probability of the path
consisting of $p_1$ and $p'$,
is at least as high as the on-time
arrival probability of the path consisting
of $p_2$ and $p'$.
Thus, a path with a travel-time distribution
that is dominated by the travel-time distribution
of at least one other path to the same node, is never
part of the optimal path from the origin $k_0$ to the destination $k^\star$.
Therefore, this subpath can be disregarded in 
subsequent iterations.
Thus, in the algorithm, there is an extra check to see if $L'$
is dominated by a distribution in the label set of $k'$,
or vice versa. Dominated
distributions are removed from the label set, and the corresponding tuples are removed from the queue.

\begin{algorithm}[ht]
 \KwResult{\doublespacing path from $k_0$ to $k^\star$ with highest on-time
 arrival probability, given departure time $t$, $B(t) = s$.}
 \textcolor{white}{a}\\
 \doublespacing
 1. With $L_0 = \{(t,1)\}$, set $D_{k_0} = \{L_0\}$
 and $D_{k_i} = \emptyset$ for all other $i$\;
 2. If $t + t_{\text{min}}(k_0,k^\star) > M$ quit
 and return 0. Else continue\;
 2. Set the queue $q = \{\{1, L_0, k_0, \text{path: }\{k_0\}\}\}$\;
 3. Extract $q^\star  = \{\alpha, L, k, P\}$ 
 with minimum first entry from $q$\;
 4. If $k = k^\star$ quit and return ($\alpha$, $P$). Else continue\;
 5. \ForEach{neighbor $k'$ of $k$ not in $P$}{
    a. Compute new travel-time distribution $L'$ via \eqref{eq:pathdist}
    with $a_{i+1}' = (k,k')$\;
    b. Compute $\alpha' = \mathbb{P}(T_{L'} + t_{\text{min}}(k',k^\star) \leqslant M)$\;
    c. \uIf{$L'$ is not dominated by an element from $D_{k'}$}
    {
    Remove
    all paths dominated by $L'$ in $D_{k'}$ and insert
    $L'$ into $D_{k'}$. Add $\{\alpha',L',k',P+k'\}$ to $q$\;
    }
 }
 6. Return to step 3.
 \caption{Optimal path for a given departure time}
 \label{alg:od1}
\end{algorithm}

\subsubsection{Optimal departure time for a path --- online setting}
We have already determined the optimal departure time $t^\star $ for a given path in the offline setting in Section~\ref{subsec:pathoffline}. Note that this departure time depends solely on the current state of the background process $B(0)$. Therefore, if case $t^\star  > 0$, the traveler should wait $t^\star $ time units until departure. However, during the time the traveler is waiting for departure, new information on the state of the background process becomes available. In the online setting, this new information is used to update the optimal departure time.
Specifically, we will consider the case where the traveler requests a new optimal departure time every $\Delta > 0$  time units. We note that the value of $\Delta$ should be chosen thoughtfully: while a higher frequency of updating the optimal departure time (i.e., a smaller value for $\Delta$) allows for more up-to-date departure time advice, the computational burden also increases.

Determining the optimal departure time using the information available at time $i\Delta$ for $i\in \mathbb{N}_0$ can then be done in a very similar way as in Algorithm~\ref{alg:offline} with some minor adjustments. Recall that the optimal departure time determined at time $i\Delta$ is written as $t^\star _{i\Delta}$. Of course, a traveler cannot indefinitely keep updating their departure time; eventually, the driver will need to depart. Whenever we find an $i$ such that
\[
t^\star _{i\Delta} - i\Delta < \Delta,
\]
it means that the optimal departure time is less than $\Delta$ time units removed from the request time. If this is the case, the driver should not wait another $\Delta$ time units to request a new departure time, since its optimal departure time is earlier than the next request time. Therefore, for such an $i$, we set $t^\star _{i\Delta}$ as the optimal departure time. 

Note that it could happen that $t^\star _{i\Delta} < i\Delta$. That is, at request time $i \Delta$, the traveler finds that they should have already departed in order to have an on-time arrival probability of at least $\eta$. In this case, the departure advice is to leave immediately at the request time $t = i \Delta$. Note that this scenario can only happen if at time $(i\!-\!1)\Delta$ it was the case that $t^\star _{(i\!-\!1)\Delta} - (i\!-\!1)\Delta \geqslant  \Delta$ as otherwise the driver should have already departed before time~$i\Delta$. Therefore, the departure advice of time~$i \Delta$ is not going to deviate more than $\Delta$ from the actual optimal departure time as defined in \eqref{optimal_departure_online}.

This procedure is now summarized in Algorithm~\ref{alg:online}.

\begin{algorithm}[ht]
 \KwResult{\doublespacing optimal departure time for traversing path $(a_1,\dots,a_m)$ within time $M$ with probability $\eta$, given $B(0) = s$.}
 \textcolor{white}{a}\\
 \doublespacing
 1. Carry out Algorithm~\ref{alg:offline}\;
 2. Set $t_0$ as the output of Algorithm~\ref{alg:offline}\;
 3. \uIf{$t_0 = 0$}{Return $t^\star  = 0$\;}
\uElse{
     Set $j = 0$\;
    \While{$t_{j\Delta}-j\Delta\geqslant  \Delta$}{
    1. Set $I_j = [\max\{0,M-j\Delta-t_{\text{max}}\}, \max\{0,M-j\Delta-t_{\text{min}}\}]$\;
    2. Use the bisection method with initial interval $I_j$ on the function
    \textbf{OnTimeProbability}($t',s, M- j\Delta$), until obtain latest departure time for which
    output is at least $\eta$. Set this departure time to $t_{j\Delta}$\;
    3. \uIf{\em{\textbf{\text{OnTimeProbability}}}$(\max\{j \Delta,M-t_{\text{max}}\},s)    \leqslant \eta$}{Return $t^\star  = j \Delta$\;}
    \uElse{
    1. Set $j = j+1$\;
    2. Obtain $B(j\Delta)$ and set $s = B(j\Delta)$\;
    }
    }
    }
    Return $t^\star  = t_{j\Delta}$
 \caption{{Optimal departure time online setting}}
  \label{alg:online}
\end{algorithm}

\section{Numerical Experiments} \label{sec:numex}
\begin{figure}[ht]
	\centering
	\includegraphics[width=0.5\textwidth]{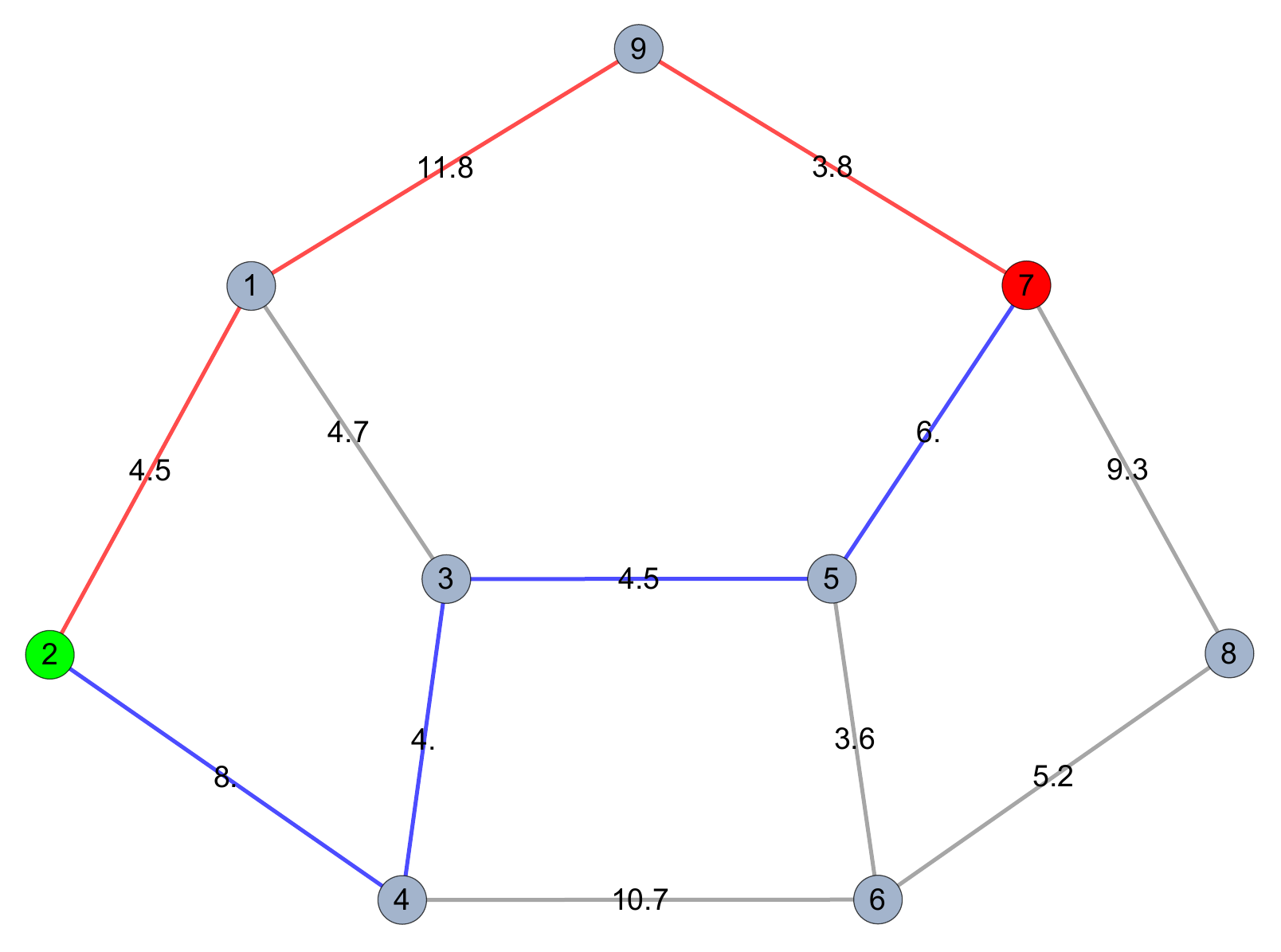}
	\caption{\doublespacing The graph used in the numerical experiments. The edge values denote the length of the edge in kilometers. We examine a driver that wants to travel from vertex 7 to vertex 2.}
	\label{fig:Amsterdam}
\end{figure}

\begin{figure}[ht]
    \centering
    \includegraphics[width=0.5\textwidth]{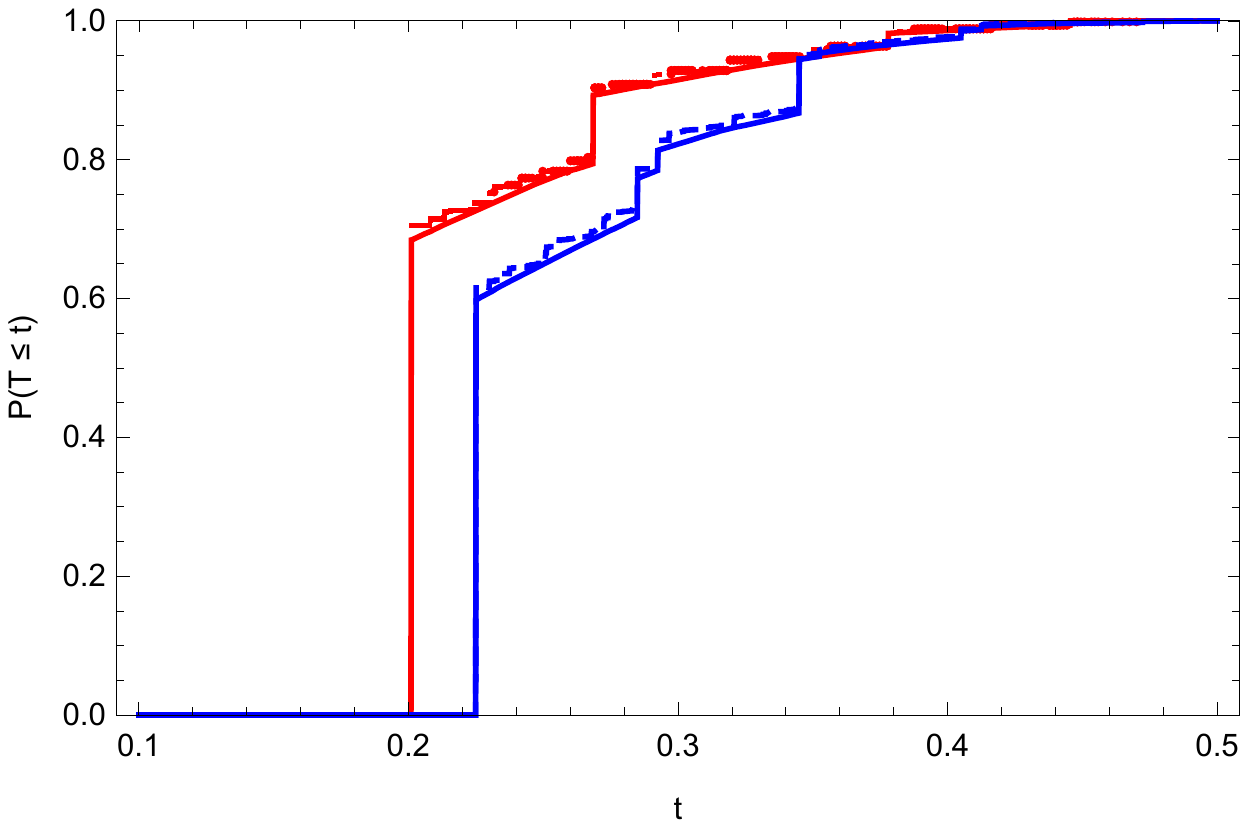}
    \caption{\doublespacing The actual (solid line) and approximated (dashed line) travel-time distribution of the red and blue route in Figure~\ref{fig:Amsterdam}.}
    \label{fig:performancediscr}
\end{figure}

Now that we have derived the optimal departure time, we will perform a set of numerical experiments in order to discuss a selection of properties of the optimal departure time, and to demonstrate the efficiency our procedure. For these experiments, we consider a road network inspired by the highways of Amsterdam, see Figure~\ref{fig:Amsterdam}.
We examine a driver that is currently located at the red-colored vertex 7 and wishes to travel to the green-colored vertex 2 by either traversing the red or blue route. We start by considering a baseline setting, which is selected for illustration purposes only and does not necessarily reflect the true parameters of this network (where we remark that experiments on a larger network  with realistic parameters are discussed in Section \ref{subsec:efficiency}). It is given by the following parameters:
\begin{itemize}
	\item[$\circ$] there are currently no incidents, which means that $B(0)=\{1\}^{12}$;
	\item[$\circ$] on each link, the velocity in case no incident occurred on that link is 100 km/h and is 40 km/h otherwise, or, $v_{a_i}(1) = 100 $ km/h and $v_{a_i}(2) = 40$ km/h for each $i=1,...,12$.
\end{itemize}
In the following experiments we may deviate from the baseline setting in order to magnify certain properties of the optimal departure time. In what follows, we note that whenever we refer to a travel-time distribution, we actually intend to refer to the \textit{approximated} travel-time distribution as described in Section~\ref{sec:optdeptime}.
For these approximations,
travel-time values are aggregated into 100 bins after
every iteration (i.e., after every step of~\eqref{eq:pathdist}). Figure~\ref{fig:performancediscr} illustrates that, for both routes in Figure~\ref{fig:Amsterdam}, the approximated travel-time distribution, obtained by the method of discretization as described in Section~\ref{sec:optdeptime}, closely resembles the actual travel-time distribution, 
obtained with 100\,000 simulation runs. {Here, on each link, both the rate of incidents and the rate of clearance are one per hour, or, $\alpha_i = \beta_i =1$ h$^{-1}$ for each $i$.}
Similar performance results were obtained for other networks
and paths under various parameter settings.
For the experiments we implemented the networks and
algorithms in Wolfram Mathematica 12.0 on 
an Intel\textregistered\ Core\texttrademark\ i7-8665U 1.90GHz computer.

\subsection{Effect of considering the on-time arrival probability} \label{subsec:expeffectprob}~\\
This experiment aims to show that, by having access to the travel-time distribution, a more suitable choice can be made in the route selection problem. {Concretely, in contrast to
only considering expected travel times, the risk-averseness of an individual driver can be incorporated, resulting in different departure times and routes. We consider the baseline setting outlined above, again with $\alpha_i = \beta_i =1$~h$^{-1}$, but with the incident rate on each link of the red route increased to two per hour.} 
Since the red route is shorter compared to the blue route, while also being more prone to incidents, it is not immediately clear which route is optimal for the driver. 

\begin{figure}[ht]
	\centering
	\includegraphics[width=0.5\textwidth]{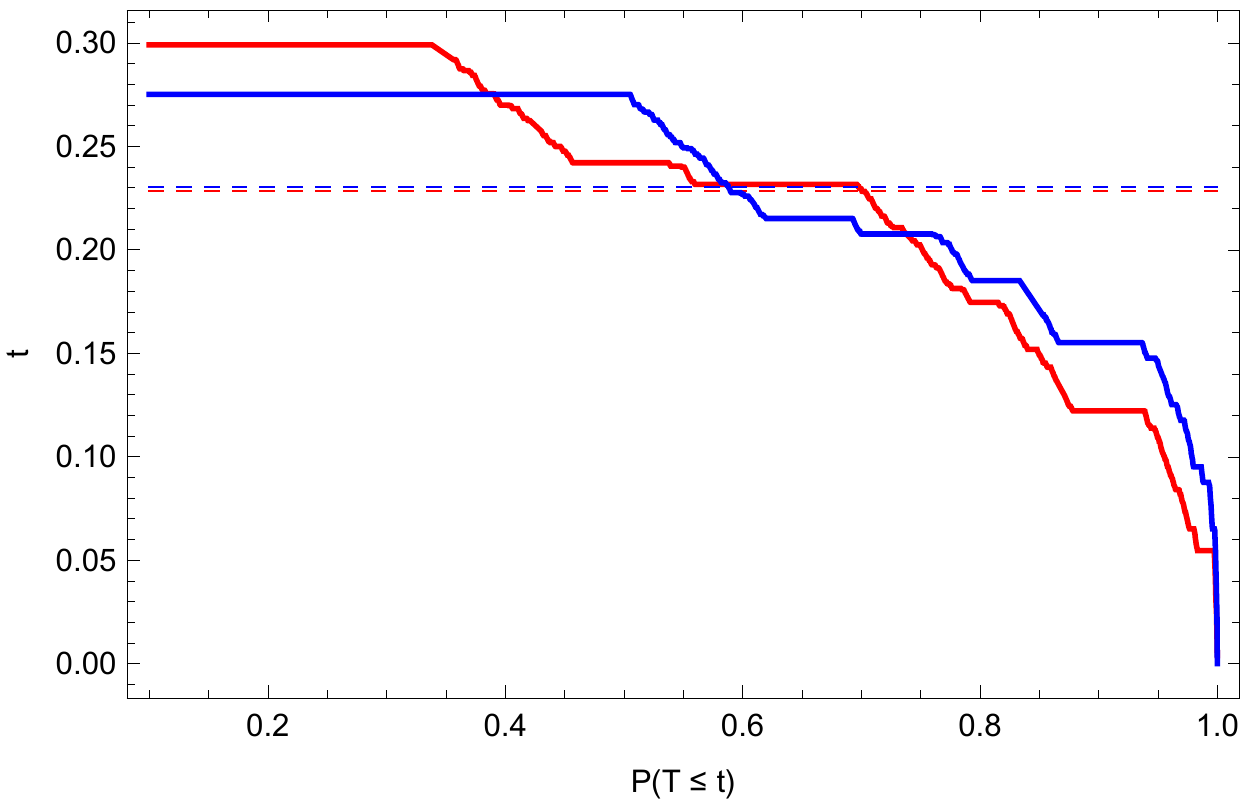}
	\caption{\doublespacing The departure time as a function of the on-time arrival probability corresponding to the red and blue route in Figure~\ref{fig:Amsterdam}, where the arrival time is chosen as $t = 0.5$. The dashed lines indicate the departure time for both routes in case of a desired arrival time at $t = 0.5$ in expectation.}
	\label{fig:on_time_arrival}
\end{figure}

Suppose that a driver wants to arrive at vertex 7 before time $t = 0.5$. We first consider the problem in which the driver wants to arrive at $t = 0.5$ in expectation. We find that the departure times for both routes in this case are very comparable, namely $t = 0.229$ for the red route and $t = 0.230$ for the blue route; see the dashed lines in Figure~\ref{fig:on_time_arrival}. In other words, if a driver were concerned with the expected arrival time, they are likely to be indifferent between the two routes.

Instead of the requirement of arriving on-time in expectation, we will now consider the problem in which a driver wants to arrive on-time at vertex 7 with a certain probability. We can do this by utilizing the arrival time distribution for both routes. From Figure~\ref{fig:on_time_arrival} we learn that the blue route allows for a later departure time, and thus should be preferred in case {the desired on-time arrival probability $\eta$ is either between 0.4 and 0.6, or greater than 0.75. In particular, for e.g. $\eta = 0.9$, the departure time of the blue route is about 0.033 h, or about 2 min, later. For routes with an approximately equal expected travel time of only 16 min, the difference between their corresponding departure times is remarkable  (i.e., it corresponds to as much as 10--15\% of the travel time).  Without having access to the travel-time distribution, which is provided by our modeling framework, this difference would have gone unnoticed.}

\subsection{The role of the background state upon departure}\label{subsec:expbackground}~\\
We consider a driver who plans to
traverse the red route in Figure~\ref{fig:Amsterdam},
with the intend to depart as late as possible
while guaranteeing an on-time arrival probability
of $\eta = 0.9$ for arrival time 14:00 hrs. We again consider the baseline setting, and we set $\alpha_i = 0.5$ h$^{-1}$ and
$\beta_i = 2$ h$^{-1}$ for all $i$. 

\begin{figure}[ht]
\centering
  \begin{subfigure}{0.47\textwidth}
  \centering
    \includegraphics[width=\textwidth]{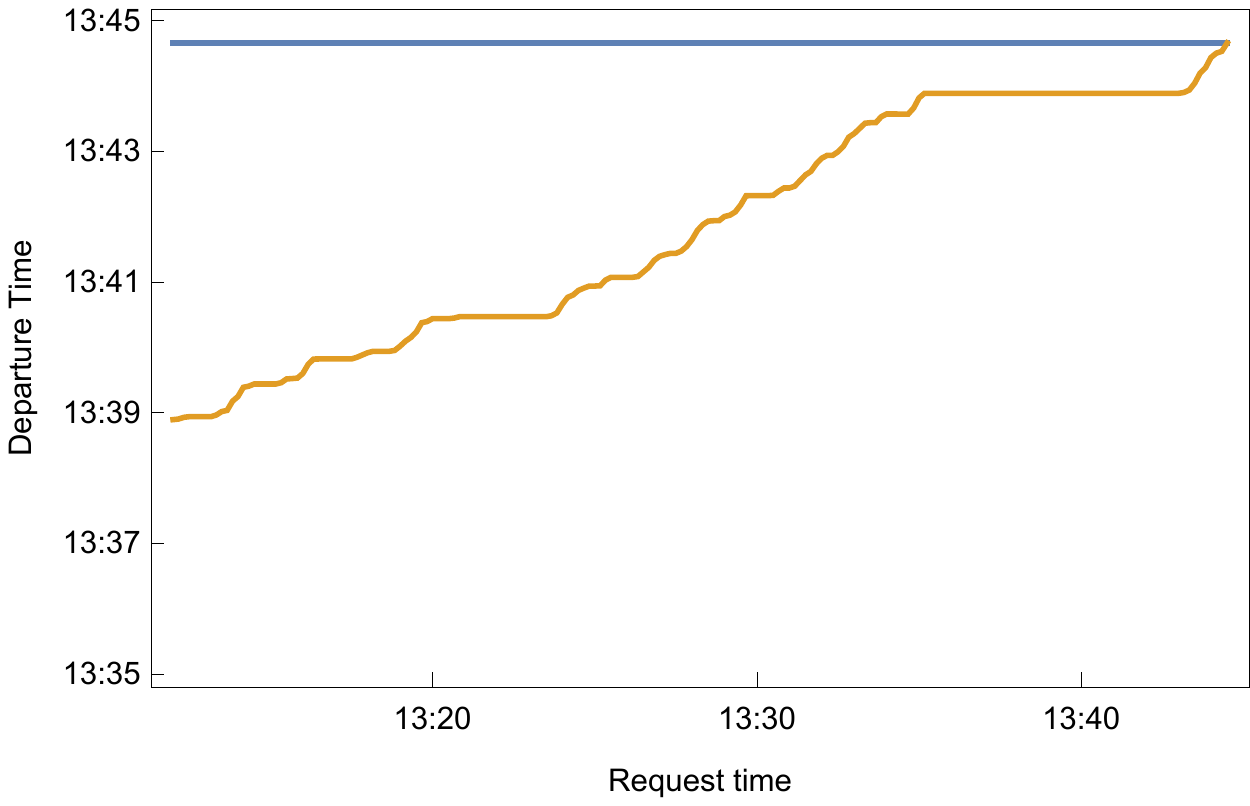}
    \caption{\doublespacing $B(0)$ is a state without incidents.}
    \label{fig:B0BtNoInc}
  \end{subfigure}
  \quad
%     \begin{subfigure}{0.1\textwidth}
%   \centering
%     \includegraphics[width=0.7\textwidth]{plotlegend.pdf}
%   \end{subfigure}
%   \quad
  \begin{subfigure}{0.47\textwidth}
  \centering
    \includegraphics[width=\textwidth]{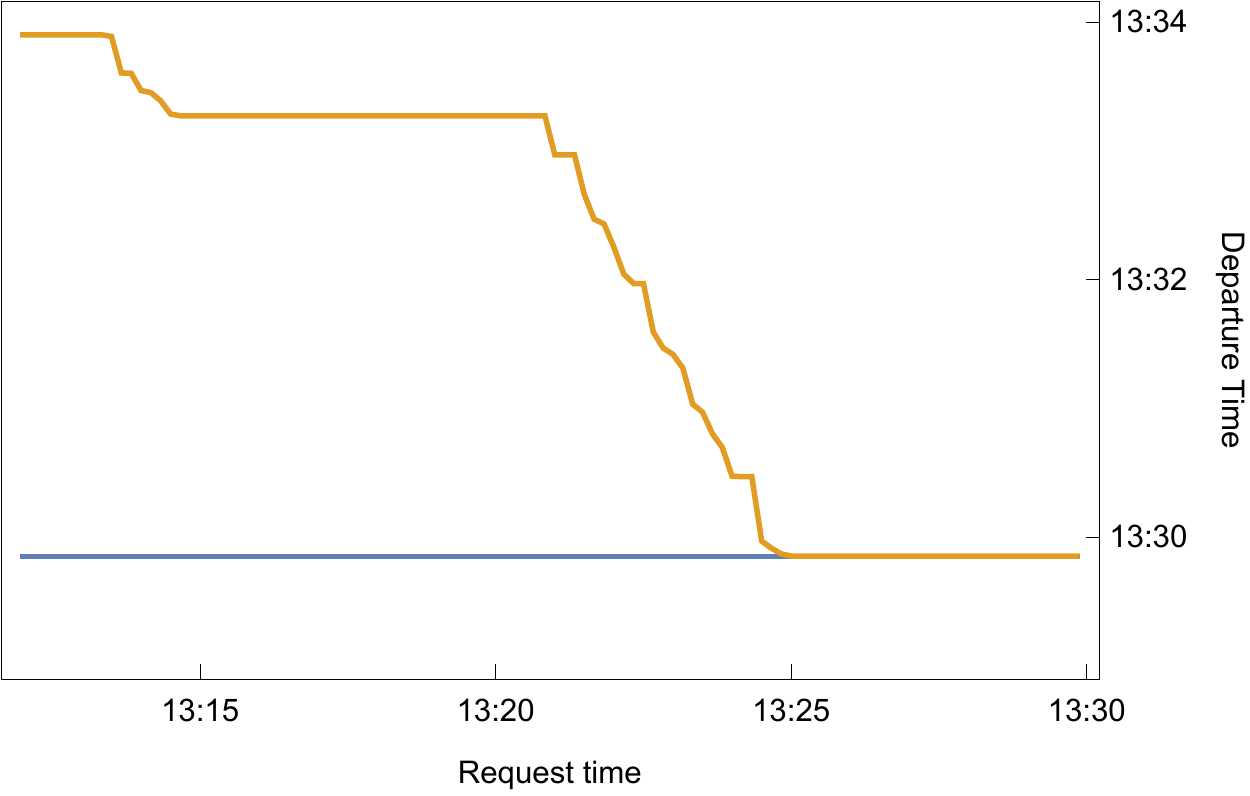}
    \caption{\doublespacing $B(0)$ is a state with solely incidents.}
    \label{fig:B0BtInc}
  \end{subfigure}
  \caption{\doublespacing Advised departure times for a vehicle that
  wants to travel the red route in Figure~\ref{fig:Amsterdam} with
  90$\%$ on-time arrival probability. Departure times are obtained
    by assuming that at departure time $t$
    the background process is in state $B(0)$ (blue),
    or by using the distribution of $B(t)$ conditional on $B(0)$ (yellow).}
  \label{fig:B0Bt}
\end{figure}

At time $0$, the moment the vehicle requests its optimal 
departure time, the background process is in a known initial 
state $B(0)$. However, for any departure time $t > 0$, the background
state $B(t)$ is random, and may very well be different from $B(0)$ due
to the occurrence or clearance of incidents in the network during the time interval~$[0,t]$. Therefore, when determining the optimal departure time $t^\star $, the possibility of  the background state having changed at time $t$ should be taken into account, which is done by using the distribution of $B(t)$ conditional on the initial state $B(0)$. The importance of doing so is illustrated in Figure~\ref{fig:B0Bt}.
% Route A: {7,9,1,2}
% Route B: {7,5,3,4,2}
% v1 = 100, v2 = 40
% alpha = 0.1, beta = 2

Figures~\ref{fig:B0BtNoInc} and~\ref{fig:B0BtInc} 
show the advised departure time of our procedure for different request times,
and compare the results with a simplified procedure, in which it is
assumed that upon departure at time $t>0$ 
the background process is still in state $B(0)$.
In the simplified procedure, the travel-time distribution is the same for every departure time $t$. Consequently, the optimal departure time under this
procedure, which we denote by \textcolor{black}{$\tilde{t}$}, is independent of the request time.
Observe that when the difference between the request time and \textcolor{black}{$\tilde{t}$} grows, 
so does the difference between 
the outputs of the two procedures. 
Specifically, in both plots, 
this difference even exceeds four minutes, 
which is substantial considering that the
travel time is in the interval of $[12.06,30.15]$ minutes.

% Notably, there is contrast between the two
% figures in what this difference
% around 13:10 implies for the functionality of 
% the simplified procedure:

In Figure~\ref{fig:B0BtNoInc}, we see that whenever the request time is before 13:44 hrs, $\tilde{t}\!>\!t^\star $.
Therefore, the on-time arrival probability for the simplified procedure will be below
the desired 90$\%$. Here,
since $B(0)$ is a state without incidents, the simplified procedure assumes that the network
is incident-free upon departure, whereas our procedure takes the possibility of changes in the background process into account. This explains why $\tilde{t}\!>\!t^\star $. The opposite phenomenon can be seen in Figure~\ref{fig:B0BtInc}. We thus see that the optimal departure time is greatly affected by the state of the background process

% In Figure~\ref{fig:B0BtInc}, around 13:10, $t_{B(0)}\!<\!t^\star $.
%     Therefore, leaving at $t_{B(0)}$ means departing earlier than necessary.
%     Note that $t_{B(0)}\!<\!t^\star $ since, with 
%     $B(0)$ a state with only incidents, the simplified procedure assumes 
%     that the network is full of incidents upon departure, whereas our procedure takes the
%     potential clearance of these incidents between request and departure into account.

% \begin{figure}[ht]
% \centering
%   \begin{subfigure}[b]{0.48\textwidth}
%     \includegraphics[width=\textwidth]{compareS0S1.pdf}
%     \caption{Departure times obtained
%     when $B(0)$ is incident-free (yellow), compared
%     to $B(0)$ with an incident on the
%     path to traverse (blue).}
%     \label{fig:compareStates}
%   \end{subfigure}
%   \quad
%   \begin{subfigure}[b]{0.48\textwidth}
%     \includegraphics[width=\textwidth]{comparePaths.pdf}
%     \caption{Departure times obtained
%     when $B(0)$ is incident-free (yellow), compared
%     to $B(0)$ with an incident on the
%     path to traverse (blue).}
%     \label{fig:compare Paths}
%   \end{subfigure}
%   \caption{Advised departure times for a vehicle that
%   wants to travel the red path of Figure~ref{[Amsterdam]} with
%   90$\%$ on-time arrival probability.}
%   \label{fig:compareStatePaths}
% \end{figure}

\subsection{Effect of request moment \& initial state} \label{subsec:exprequestmoment}
\begin{figure}[ht]
  \begin{center}
    \includegraphics[width=0.65\textwidth]{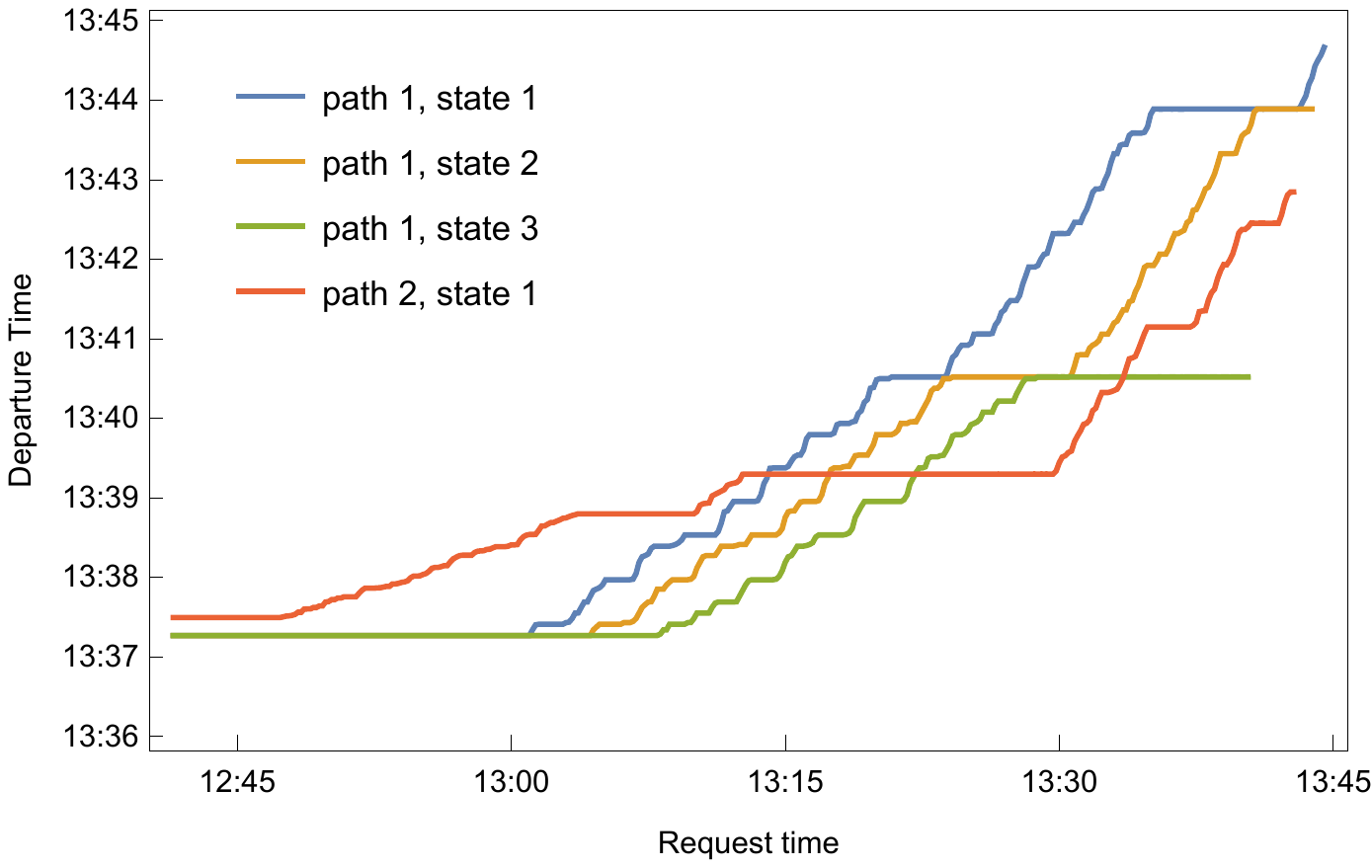}
  \end{center}
  \caption{\doublespacing Departure times for vehicle that
    wants to travel the red path (path~1) or
    the blue path (path~2) in Figure~\ref{fig:Amsterdam}, and wants to
    arrive before 14:00h with 90$\%$ certainty. We consider
    three initial states: incident-free (state~1), 
    only an incident at link $(1,2)$ (state~2)
    and only an incident at link $(7,9)$ (state~3).}
  \label{fig:comparePaths}
\end{figure}

~\\
% Plek van het ongeluk is belangrijk: Dit lijkt op
% het vorige experiment. Laat zien dat de toestand van het netwerk in
% het begin van de route een grotere rol speelt in het selectieproces.
% In steady state: 1 optimale route
% Maar als dichter erop: kan ook andere route
Besides showing the importance of incorporating the dynamics of $B(t)$,
Figures~\ref{fig:B0BtNoInc} and~\ref{fig:B0BtInc} also display that $t^\star $ depends
on the initial state upon the time of request.
This effect is even more clear from
Figure~\ref{fig:comparePaths}, which shows the optimal 
departure time $t^\star $ for traversing the red path (path~1)
of Figure~\ref{fig:Amsterdam} under three different initial states. Note that, similar
to the preceding experiment, the driver wishes to arrive
before 14:00 hrs with 90$\%$ certainty. Moreover, we again
let $\alpha_i = 0.5$ h$^{-1}$ and $\beta_i = 2$ h$^{-1}$.

First, we note that, being in steady state, a request time before 12:45 hrs yields similar departure times for
the three initial states. This is no longer the case when the request moment is
close to the arrival time.
As expected, the latest departure time is obtained for the
incident-free background state. Importantly, in case there is an incident
at request time, the location of this incident has a significant impact on
the departure time $t^\star$.
That is, with the location of the incident at the end
of the path (state~2), there is a high probability of
clearance before arrival at the incident link, such that
the corresponding departure time 
is still relatively close to the departure time of
the incident-free state. However, in case the location
of the incident is at the start of the path (state~3),
there is only a small probability that this incident
is cleared upon arrival at the incident-link. Consequently,
the vehicle must depart considerably earlier in state~3 than in state~2.

% \begin{wrapfigure}{r}{0.60\textwidth}
%   \begin{center}
%     \includegraphics[width=0.95\textwidth]{comparePathsV2.pdf}
%   \end{center}
%   \caption{\doublespacing Departure times for vehicle that
%     wants to travel the red path (path~1) or
%     the blue path (path~2) in Figure~\ref{fig:Amsterdam}, and wants to
%     arrive before 14:00h with 90$\%$ certainty. We consider
%     three initial states: incident-free (state~1), 
%     only an incident at link $(1,2)$ (state~2)
%     and only an incident at link $(7,9)$ (state~3).}
%   \label{fig:comparePaths}
% \end{wrapfigure}
In case the driver would only give the OD-pair, there
are other paths from node 7 to node 2
that potentially outperform the red path,
of which the the blue path of Figure~\ref{fig:Amsterdam} (path~2)
is an example.
% Naturally, if there is an incident on path~1, this driver
% should also consider other paths that start at node 7
% and have node 2 as destination. An example of such a path
% is the blue path of Figure~\ref{fig:Amsterdam} (path~2). 
% First, 
Note that as the incidents in state~2 and~3 are not located on path~2,
all three initial states considered in Figure~\ref{fig:comparePaths} 
will result in the same departure time for traversal via path~2.
% The total distance of path~2 is slightly longer, but, 
As is shown in
Figure~\ref{fig:comparePaths}, the vehicle
would prefer path~2 over path~1 in steady state.
% This can be explained
% by the fact that path~2 consists of more links than path~1,
% such that an incident on one of those links has less effect
% on the total travel time of the path. 
% However, this dominance
% is no longer visible when the request time is close
% to the arrival time, i.e., when the probability of occurrence
% of incidents is small, and consequently, their effect
% on the travel times on the path. Then, path~2 is only dominant
% over path~1 if there is an incident located at the start
% of path~2. In case path~2 is free of incidents, or has
% an incident at the final link, the optimal departure time
% is further in the future than the optimal departure time
% under path~2.
However, this dominance
no longer applies when the request time is close
to the arrival time, as, in that case, path~2 is only preferred
over path~1 if there is an incident located at the start
of path~1.

{From this experiment we learn that our modeling procedure effectively exploits knowledge of the locations of the currently present non-recurrent events in the network. It not only incorporates the presence of incidents on the routes, but also distinguishes between the locations of these incidents. For example, incidents on a link at the end of a path affect the departure time to a lesser extent, as our framework incorporates the high probability of incident clearance before arrival at this link. This is of course not the case for incidents closer to the origin.}

\subsection{Online vs Offline Departure Time} \label{subsec:onlinevsoffline}~\\
{As our modeling framework allows for the real-time implementation of the changing conditions in the road
network, we now implement the online version of the optimal departure time problem outlined in Algorithm~\ref{alg:online}, in which the traveler receives departure time updates while still at the origin.} We again consider the baseline setting defined in the beginning of Section~\ref{sec:numex}, with the rate of incidents $\alpha_i = 0.1$ h$^{-1}$ and the rate of clearances $\beta_i = 2$ h$^{-1}$. We will study a driver that wants to travel from vertex 7 to vertex 2 using the red route in Figure~\ref{fig:Amsterdam}.

\begin{figure}[ht]
    \centering
    \includegraphics{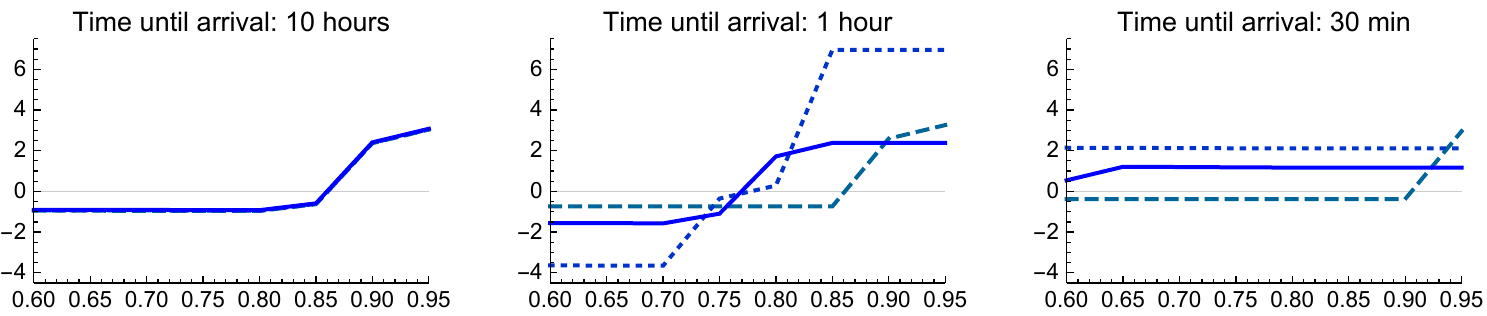}
    \caption{\doublespacing The difference between the expected online departure time, approximated using 10.000 repetitions, and the offline departure time. This is done for a selection of arrival times $M$ and background states $B(0)$, where the dashed line corresponds to no incidents, the solid line to an incident on link $(1,2)$, and the dotted line to an incident on link $(7,9)$.}
    \label{fig:online}
\end{figure}

We compute the optimal departure of this driver using both the online and the offline setting, utilizing Algorithm~\ref{alg:online} and Algorithm~\ref{alg:offline}, respectively. We do this for a selection of requested arrival times: we let the length of the interval between request and desired arrival time be $M = 10$ h, $M = 1$ h, and $M = 30$ min. In addition, we consider a range of on-time arrival probabilities $\eta$. As the online optimal departure time is random, we approximate its expectation by performing $10\,000$ repetitions and computing its mean. Lastly, we subtract the offline departure time from the approximated expected online departure time. The findings are shown in Figure~\ref{fig:online}.

In the left graph of Figure~\ref{fig:online} (that corresponds to $M= 10$ h), we see that the online setting gives a later departure time compared to the offline setting in case the driver demands a high on-time arrival probability ($\eta> 0.85$). This difference in departure time of 3 minutes is quite substantial, as the expected travel time of this route is roughly 15 minutes (i.e., 20\% of the duration of the trip). This can be explained as follows. Since the requested arrival time is 10 hours from now, and since only the current state of the network is known in the offline setting, the actual state of the network on departure is highly uncertain. Therefore, if a driver demands a high on-time arrival probability in the offline setting, a very conservative believe about the state of the network on time of departure must be used. Either this belief is not far off, and the online and offline departure time will not much differ, or this belief was indeed too pessimistic and the network is actually in a more favorable state upon departure such that the online setting will give a later optimal departure time. On average, we see that the departure time corresponding to the online setting will be later than that of the offline setting. Note that the converse is also true: for lower values of the on-time arrival probability, the online setting will find an earlier departure time compared to the offline setting. 
%The case where a low on-time arrival probability is requested, however, seems of less relevance. 

With this in mind, it is not surprising that the difference between the online and offline optimal departure time is smaller when $M =1$ h, and is even smaller in the $M = 30$ min case. As the time of request moves closer to the departure time, the offline setting is able to determine the optimal departure time based on a more recent state of the network. This way, the state of the network upon departure is less uncertain and the offline method will more closely resemble the online method.

Lastly, in Figure~\ref{fig:online}, we also consider the effect of the initial background state $B(0)$ on the difference between the optimal departure times for the online and offline setting. For $M=10$ h, the initial background state plays no role, as the distribution of the background state upon departure, given any initial state, is effectively in the steady state. This is not the case for $M=1$ h and $M=30$ min, and we see that the difference between the online and offline procedure is greatest for the case in which there is an incident on the nearby link $(7,9)$ (the dotted line), and smallest in case there are no incidents (the dashed line).

As expected, more is to be gained from the online procedure when there is more uncertainty about the state of the network upon departure, which is in particular the case when the initial state contains incidents. Especially when these incidents are located near the origin, such as link $(7,9)$, the online method is able to more accurately incorporate the state of nearby links as they will be seen by the traveler upon departure, as opposed to links that are further away from the origin, such as link $(1,2)$.

% \begin{table}[ht]
%     \centering
%     \begin{tabular}{|c||c|c|c|}
%     \hline
%      & 30 min & 1 h & 10 h \\ \hline
%     State 1 & 100$\%$ & 96.9$\%$ & 96.9$\%$\\
%     State 2 & 100$\%$ & 97.1$\%$ & 97.1$\%$\\
%     State 3 & 100$\%$ & 97.0$\%$ & 97.0$\%$\\ \hline
%     \end{tabular}
%     \caption{10 000 experiments, per setting the minimum
%     success rate among the different considered values
%     of $\eta$}
%     \label{tab:my_label}
% \end{table}

%dashed line: no accidents
%dotted line: accident on link 2
%solid line: accident on link 3

\begin{figure}[ht]
    \centering
    \includegraphics[width=0.4\textwidth]{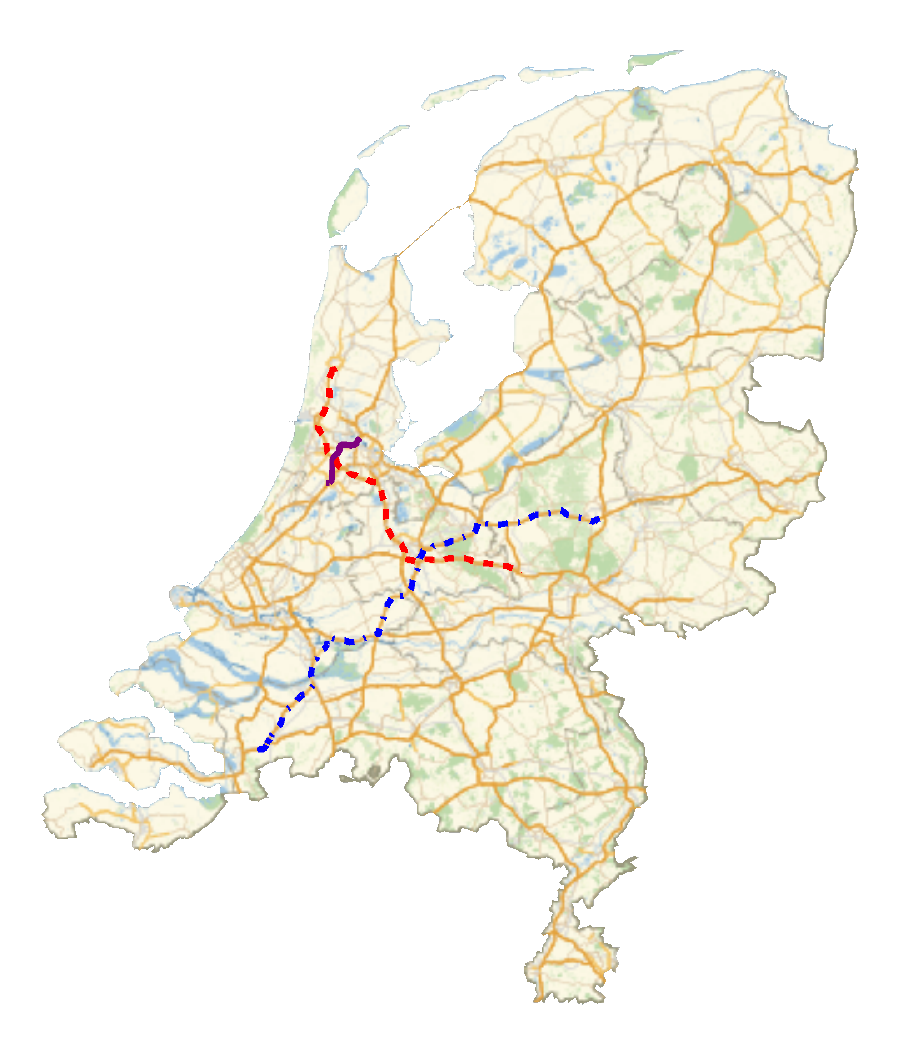}
    \caption{\doublespacing The Dutch highway system, with blue, purple and red path highlighted.}
    \label{fig:longpath}
\end{figure}

\subsection{Efficiency} \label{subsec:efficiency}~\\
To empirically study the computational complexity of our procedures, we 
consider the complete Dutch highway network as depicted
in Figure~\ref{fig:longpath}. This network consists of 659
nodes (i.e., highway ramps)
and 1378 links (i.e., roads connecting ramps).
The parameters are chosen to reflect the real-life setting
and estimated via the fitting procedures presented
in 
\cite{Levering2022Estimation}.

The computational performance of Algorithm~\ref{alg:offline}
is given in Table~\ref{tab:comptimeofflinepath}. 
To this end, we consider
the blue path in Figure~\ref{fig:longpath}, 
%from Apeldoorn
%to Roosendaal, 
and compute the optimal departure time 
for subpaths of different lengths from the origin.
It can be observed that the run-time of the algorithm
increases super-linearly in the number of links in the path.

\begin{table}[b]
    \centering
    \begin{tabular}{|c|c|c|c|c|c|}
    \hline
    Number of links & 1 & 5 & 10 & 15 & 20 \\ \hline
    Seconds & $<10^{-4}$ & 0.46 & 1.36 & 4.05 & 6.60 \\
    \hline
    \end{tabular}
    \caption{\doublespacing Computational performance optimal
    departure time path (in sec.).}
    \label{tab:comptimeofflinepath}
\end{table}

% \begin{table}[ht]
%     \centering
%     \begin{tabular}{|c|c|c|c|c|c|}
%     \hline
%     OD-pair distance & 5 & 10 & 15 & 25 & 50 \\ \hline
%     Algorithm~\ref{alg:od1} & & & & & \\
%     \hline
%     $k$-shortest paths & & & & &  \\
%     \hline
%     \end{tabular}
%     \caption{Computational performance optimal
%     departure time OD-pair (in sec.)}
%     \label{tab:comptimeofflineodpair}
% \end{table}

For the performance of the two procedures
proposed in Section~\ref{subsub:od}, designed to 
output the optimal path and departure
time in case a traveler gives the
desired arrival time for an OD-pair, we consider
% To this end, we consider various OD-pairs in the 
% Dutch highway network. The OD-pairs differ in
% length, which we measure by the minimum number of links
% a vehicle needs to traverse to get from the origin
% to the destination. For every OD-pair distance, we
% randomly pick twenty OD-pairs from the network, and compute the 
% mean run-time of a procedure over these twenty pairs. 
two OD pairs that differ considerably in length (i.e.,
minimum number of links a vehicle needs to traverse to
get from the origin to the destination). 
Concretely, OD-pair 1 has a length of 5 
(Figure~\ref{fig:longpath}, purple path), whereas OD-pair
2 has a length of 32 (Figure~\ref{fig:longpath}, red path).
For the experiments, we look at a traveler with an on-time arrival 
probability objective of $0.8$, and set
$M = 1$ h  for OD-pair 1 and and $M = 2$ h for OD-pair 2. 
Now, as only incidents located in an area around 
(one of) the shortest paths
between an OD-pair will significantly affect the departure time 
and route advice, we focus on the states $B(0)$ that 
encode incident instances on one of the three shortest
paths (in km) between the origin and destination.
With the probability of more than two incident
instances being low, we further restrict our focus
and only consider the set of states $B(0)$ that reflect one
or two incident occurrences on the set of three shortest paths
between an origin and destination, and use $B^\star $
to denote this set. The setting
in which the three shortest paths are
completely incident-free is treated separately.

Recall that, contrary
to the $k$-shortest path method, 
the procedure which employs
bisection on the label correcting
algorithm from Algorithm~\ref{alg:od1} is exact,
and outputs the optimal departure time and corresponding
path. Now, for both OD pairs, 
whenever $B(0)$ is such that there is no incident
on one of the three shortest paths, the $k$-shortest path
method already yields the optimal
path and departure time for $k=1$.
Note that this setting will be frequently encountered: the stationary probabilities of such initial
states are $0.994$ and $0.985$ for OD pairs~1 and~2, respectively.
In case $B(0) \in B^\star $,
%is such that there are incident(s) on one of
%the three shortest paths, 
$k=1$ is typically not sufficient for
obtaining the optimal departure time. 
To assess the performance for different values of $k$,
we compute, for every state $B(0) \in B^\star $ 
%that encodes
%the occurrence of one or two incidents on
%the set of three shortest paths between an OD-pair, 
the
optimal departure time. A quantification of the
performance is expressed in terms of the mean absolute percentage error (MAPE), defined
as
\begin{align*}
    \text{MAPE}: = \frac{1}{|B^\star |}\sum_{s \in B^\star } 
    \frac{|t_s^\star -t_s^k|}{t_s^\star },
\end{align*}
with $t_s^\star $ the optimal departure time and
$t^k_s$ the departure time as outputted by the $k$ shortest path
procedure, given $B(0) = s$.

\begin{table}[t]
    \centering
    \begin{tabular}{|c|l||c|c||c|c|}
    \hline
    \multicolumn{2}{|c||}{} & \multicolumn{2}{c||}{Run-time (sec.)} & \multicolumn{2}{c|}{MAPE} \\
    \hline
    \multicolumn{2}{|c||}{OD-pair } & 1 & 2  
    & 1 & 2 \\ 
    \hline
    \multicolumn{2}{|c||}{Algorithm~\ref{alg:od1}} & 4.0 & 393.8
    & \multicolumn{2}{c|}{0$\%$}\\
    \hline
    \multirow{5}{*}{$k$-shortest paths} 
    & $k=1$ & 0.9 & 3.1
    & 4.5$\%$ & $2.9\%$ \\
    & $k=3$ & 1.9 & 9.6
    & 0.0$\%$ & $0.7\%$ \\
    & $k=5$ & 4.8 & 16.2
    & 0.0$\%$ & $0.5\%$ \\
    & $k=10$ & 16.2 & 34.6
    & 0.0$\%$ & $0.4\%$ \\
    & $k=15$ & 26.0 & 55.7
    & 0.0$\%$ & $0.0\%$ \\ \hline
    \end{tabular}
    \caption{\doublespacing Computational performance optimal
    departure time OD-pair (in sec.).}
    \label{tab:comptimeofflineodpair}
\end{table}

The procedure based on Algorithm~\ref{alg:od1} is exact, and thus has
zero MAPE. The procedure based on the $k$-shortest
path algorithm is not exact, but, as can be observed from
Table~\ref{tab:comptimeofflineodpair}, already
provides near-optimal results for $k = 3$. 
Moreover, the computational savings when using
the $k$-shortest path method with moderate values of $k$
is significant, especially for OD-pair~2.
{
Note that, in these experiments, we have only used
one core, whereas the run-time of the $k$ shortest path
algorithm can be further reduced when applying
parallel computing, as the optimal departure times
of the $k$ paths are computed individually.
% parallel comment
% Lastly, we consider an alternative
% objective function, 
% where 
% -- in contrast
% to the moment of request --
% a vehicle
% wants an on-time arrival probability guarantee
% at the moment of departure. We show that this
% objective function can be computed
% using the techniques presented in
% Section~\ref{sec:optdeptime}, and compare
% the output and speed of the method with 
% the computation of the optimal departure
% time as defined in Section~\ref{sec:problem}.
Importantly, in Appendix~\ref{app:speedup}, we argue
that, in case of the $k$ shortest path method, 
even more computational efficiency can be realized
by allowing additional precomputations.
}

\section{Conclusive remarks} \label{sec:concl}
In this paper we  have developed a procedure for determining the optimal departure time in road networks with stochastic disruptions. We accomplished this by first developing an iterative procedure by which the travel-time distribution can be numerically evaluated for each departure time. Then, these distributions enabled us to develop efficient algorithms that identify the optimal departure time. We performed a selection of numerical experiments that exemplify various properties of the optimal departure time, and we demonstrated the efficiency of our procedure by applying it to an existing (large) road network. Lastly, we outlined an extension of our framework to traffic networks with more general dynamics, and we provided a speed-up technique for our procedure.

We defined the optimal departure time as the latest  time of departure such that a selected on-time arrival probability is still guaranteed. By allowing the departure time to depend on the on-time arrival probability, the risk averseness of drivers can be taken into consideration. Next, in order to account for both recurrent and non-recurrent congestion, we used the Markovian Velocity Model that relies on a background process that tracks the events affecting the velocities in the traffic network. Doing so, our model also successfully exploits the knowledge of the locations of the currently present events in the road network. This allowed us to develop an online version of the optimal departure time problem, in which the traveler (while still at the origin) receives departure time updates, that incorporate the most recent state of the road network.

Our numerical experiments illustrated the extent to which the optimal departure time is affected by the state of the background process and the time of request. We also demonstrated that the route selection process is affected by using the latest departure time as an objective function. Moreover, we were able to quantify the substantial reduction in travel time budget that can be obtained by utilizing the online version of the problem. Lastly, we have demonstrated that our procedure can also be successfully employed in a real-world road network, as the run-time of our procedure, even in large road networks, remains manageable.

Several directions for follow-up research can be thought of. First, one could focus on empirically validating the approach presented in this paper. Secondly, one may develop an interface by which a driver can, explicitly or implicitly, reveal their risk averseness, after which the optimal departure time can be communicated to the driver. Thirdly, our framework could be extended to a setting that also allows for adaptive routing once the driver has already departed. Lastly, different notions of optimality for the departure time, depending on the travel-time distribution, may be considered. An example could be a variant that, in case of late arrival, also takes into account by how much the desired arrival time has been exceeded.

\section*{Acknowledgements}
 The authors would like to thank dr.\ Marko Boon (TU Eindhoven)
for his input on the implementation of the procedures,
as well as the helpful feedback on the manuscript.

% The authors would like to thank dr.\ Maaike Snelder (TU Delft)
% for her helpful feedback on the manuscript as well as her suggestions
% on the analysis of the NDW data.

\bibliographystyle{agsm}
\bibliography{references2.bib}

\appendix
% The model and methods we have presented
% allow for various extensions, of which we will
% discuss two types. First, 
% we provide examples to show how the algorithms of 
% Section~\ref{sec:optdeptime} can be adapted
% in case the {\sc mvm} is more comprehensive,
% e.g.\ allowing for recurrent events
% or non-exponential incident durations.
% Second, we focus on the
% speed of the algorithms, and
% present a way to reduce the run-time of 
% Algorithm~\ref{alg:offline} and 
% the $k$-shortest path method
% in case the paths consist of a high number of links,
% thus making these procedures real-time.

\section{Recurrent congestion and non-exponential incident duration}\label{app:recurrent}~\\
For expositional reasons we have considered a
compact version of the {\sc mvm}. That is,
the algorithms were designed under the setting
that the only events affecting arc speeds are
incidents and, moreover, that both the
duration of these incidents and the time between
these incidents are exponentially distributed.
In Remark~1, we have claimed that we can handle
a more comprehensive version of the {\sc mvm} as well.
To corroborate this claim, we will provide two examples that
show how to adapt the algorithms of 
Section~\ref{sec:optdeptime} in case a more detailed
version of the {\sc mvm} is used. 
In the first example, we consider
a setting in which the per-link incident duration
is not exponentially distributed, whereas
in the second example, 
we discuss including recurrent
events into the departure time framework.

\begin{example} \label{ex:nonexponential}
{\em {
In the current model setting, the duration of
an incident follows an exponential distribution,
although generally, this is not always a realistic assumption.
In \cite{Levering2022Estimation},
in which the {\sc mvm} is used to
describe travel times in the Dutch highway network,
it was e.g.\ shown that there are indeed links for which
the exponential distribution does not represent 
the incident duration well. It is moreover argued,
that, for these links,
the incident duration is well described 
by an Erlang-2, mixture-Erlang or
hyperexponential distribution.
{
Below we argue how to obtain the optimal
departure time in case the incident duration
is indeed modeled by one of these distributions.
% These being phase-type distributions, they
% fit into the {\sc mvm} framework. 
}

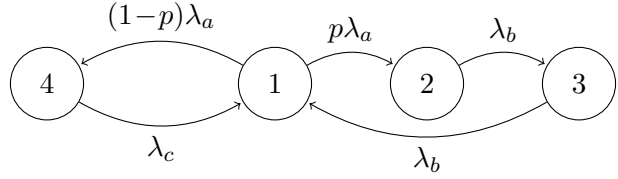
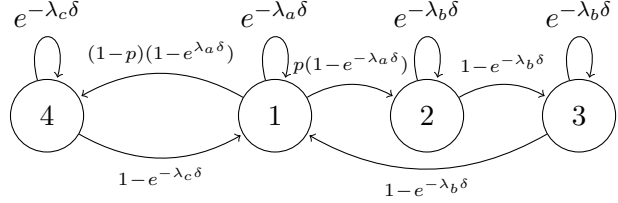
\begin{figure}[ht]
\centering
  \begin{subfigure}[b]{0.55\textwidth}
\begin{centering}
\begin{tikzpicture}[shorten >=1pt,node distance=2cm,on grid,auto]
    \node[state] (q_0) {$1$};
    \node[state] (q_1) [right=of q_0] {$2$};
    \node[state] (q_2) [right=of q_1] {$3$};
    \node[state] (q_3) [left= 3 of q_0] {$4$};

    \path[->]
    (q_0) edge [bend left] node {$p\lambda_{a}$} (q_1)
    (q_1) edge [bend left] node {$\lambda_{b}$} (q_2)
    (q_2) edge [bend left] node {$\lambda_{b}$} (q_0)
    (q_0) edge [bend right] node [above] {$(1\!-\!p)\lambda_{a}$} (q_3)
    (q_3) edge [bend right] node [below] {$\lambda_{c}$} (q_0);
\end{tikzpicture}
    \caption{Continuous.}
    \label{fig:mixtureerlangCont}
      \end{centering}
  \end{subfigure}
  \\
  \begin{subfigure}[b]{0.55\textwidth}
  \begin{centering}
\begin{tikzpicture}[shorten >=1pt,node distance=2cm,on grid,auto]
    \node[state] (q_0a) {$1$};
    \node[state] (q_1a) [right=of q_0a] {$2$};
    \node[state] (q_2a) [right=of q_1a] {$3$};
    \node[state] (q_3a) [left= 3 of q_0a] {$4$};

    \path[->]
    (q_3a) edge [loop above] node {$e^{-\lambda_c \delta}$} (q_3a)
    (q_2a) edge [loop above] node {$e^{-\lambda_b \delta}$} (q_2a)
    (q_1a) edge [loop above] node {$e^{-\lambda_b \delta}$} (q_1a)
    (q_0a) edge [loop above] node {$e^{-\lambda_a \delta}$} (q_0a)
    (q_0a) edge [bend left] node {{\tiny $p(1\!-\!e^{-\lambda_{a}\delta})$}} (q_1a)
    (q_1a) edge [bend left] node {{\tiny $1\!-\!e^{-\lambda_{b}\delta}$}} (q_2a)
    (q_2a) edge [bend left] node {{\tiny $1\!-\!e^{-\lambda_{b}\delta}$}} (q_0a)
    (q_0a) edge [bend right] node [above] {{\tiny $(1\!-\!p)(1\!-\!e^{\lambda_{a}\delta})$}} (q_3a)
    (q_3a) edge [bend right] node [below] {{\tiny $1\!-\!e^{-\lambda_{c}\delta}$}} (q_0a);
\end{tikzpicture}
    \caption{Discrete.}
    \label{fig:mixtureerlangDisc}
\end{centering}
  \end{subfigure}
    \caption{\doublespacing Structure of $X_a(t)$ whenever the incident
    duration has a mixture Erlang-2 (with probability~$p$) 
    and Erlang-1 (with probability~$1-p$) distribution,
    with $\lambda_a, \lambda_b, \lambda_c \in \mathbb{R}_{> 0}$
    and $p \in (0,1)$. State 1
    represents the incident free phase.}
    \label{fig:mixtureerlang}
\end{figure}

In this example, we consider the setting
in which the incident distribution on a link $a \in A$
is a mixture of an Erlang-2 (with probability $p$) 
and Erlang-1 (with probability $1\!-\!p$) distribution,
as depicted in Figure~\ref{fig:mixtureerlangCont}.
Now, given $X_a(0) = j$, the link
travel-time distribution of $a$ can be found
with a procedure in a similar
spirit as Algorithm~\ref{alg:distrsinglelink},
using the discrete Markov chain in
Figure~\ref{fig:mixtureerlangDisc}.
% Indeed, for a given (state, distance, probability)-tuple, 
% it is known what the potential
% states in the next time step are, with
% which probability they arise, and what distance
% one travels when driving with the 
% speed corresponding
% to these states.
Let $L_a^j$ again denote the list
of (travel time, probability)-pairs that
arises.
{Note that, as there are now five potential
non-symmetric
transitions instead of two symmetric ones,
in determining the granularity, 
$\delta$ should satisfy five conditions
of the form in \eqref{eq:granerrorapprox1},
as well as five conditions of the form in~\eqref{eq:granerrorapprox2}.
%are now maxima over five terms.
}

In a similar fashion, we can determine approximate
link travel-time distributions for 
{any mixture Erlang, \mbox{Erlang-2} or
hyperexponential incident duration.}
% the Erlang-2 and
% hyperexponential incident durations. 
Then, computation
of the path travel time is straightforward, as we 
can use \eqref{eq:pathdist} to compute the travel time
on a path, replacing $j=1,2$ by an iteration over 
the number of states in the state space of the link.
A similar adaptation allows the use of 
Algorithm~\ref{alg:offline} for determining the optimal departure
time of a path, and the use of both proposed methods
for the optimal departure time of an OD-pair.

% The only difficulty one encounters when working with
% phase-type incident durations is the identification
% of $B(0)$. That is, if there is an incident upon departure,
% it is not immediately known in which phase the incident 
% is. However, knowledge of the current running time of the
% incident allows the computation of the distribution 
% vector of $B(0)$, which will then replace
% $\boldsymbol{p}_0^s$ in \eqref{distBOffline}.
% Indeed, for the process depicted in Figure~\ref{fig:mixtureerlang}, 
% in case the incident duration $Y$ has survived
% a time $t > 0$, we have for $j = 2,3,4$:
% \begin{align*}
%     \mathbb{P}(X_a(0) = j | Y > t, X_a(0) \neq 1) = 
%     \frac{\mathbb{P}(Y > t | X_a(0) = j)\mathbb{P}(X_a(0) = j | X_a(0) \neq 1)}{
%     \sum_{i = 2}^4 \mathbb{P}(Y > t | X_a(0) = i)\mathbb{P}(X_a(0) = i | X_a(0) \neq 1)}.
% \end{align*}
% This fraction can be computed by realizing that
% for $E_1, E_2 \sim \text{Exponential}(\lambda_2)$
% and $E_3 \sim \text{Exponential}(\lambda_3)$
% \begin{align*}
%     \mathbb{P}(Y > t | X_a(0) = 2) &= \mathbb{P}(E_1 > t) = e^{-\lambda_2t} \\
%     \mathbb{P}(Y > t | X_a(0) = 3) &= \mathbb{P}(E_1 + E_2 > t) = (1+\lambda_2t)e^{-\lambda_2t}\\
%     \mathbb{P}(Y > t | X_a(0) = 4) &= \mathbb{P}(E_3 > t) = e^{-\lambda_3t}
% \end{align*}
% and
% \begin{align*}
%     \mathbb{P}(X_a(0) = i | X_a(0) \neq 1) &= p/2 
%     \quad \quad \quad \quad i = 2,3\\
%     \mathbb{P}(X_a(0) = 4 | X_a(0) \neq 1) &= 1-p.
% \end{align*}

The only difficulty one encounters when working with these more general 
phase-type incident durations is the identification
of $B(0)$. The complication is that, if there is an incident upon departure, the current state of the network is unknown as there are multiple states that correspond to the occurrence of an incident. However, knowledge of the elapsed incident duration  allows the computation of the distribution
vector of $B(0)$, which will then replace
$\boldsymbol{p}_0^s$ in \eqref{distBOffline}.
Indeed, writing $Y$ for the duration of the current incident,
$t$ for the (known) elapsed duration of the incident, and $F_1$ ($F_2$, respectively) 
for the event that the incident duration has an 
Erlang-$1$ (Erlang-$2$, respectively) distribution,
we have, for $j = 2,3$:
\begin{align*}
    \mathbb{P}(X_a(0) = j \,|\, Y \geqslant  t) &= \mathbb{P}(X_a(0) = j \,|\, Y \geqslant  t, F_2)\mathbb{P}(F_2 \,|\, Y \geqslant  t) \\
    &= \frac{p\,\mathbb{P}(X_a(0) = j \,|\, Y \geqslant  t, F_2)\,\mathbb{P}(Y \geqslant  t \,|\, F_2)}{
    (1-p)\,\mathbb{P}(Y \geqslant  t \,|\, F_1) + p\,\mathbb{P}(Y \geqslant  t \,|\, F_2)}.
\intertext{This fraction can be computed by realizing that
for $E_1, E_2 \sim \text{Exponential}(\lambda_b)$
and 
\protect\linebreak $E_3 \sim \text{Exponential}(\lambda_c)$, we have}
    \mathbb{P}(Y \geqslant  t \,|\, F_1) &= \mathbb{P}(E_3 > t) = e^{-\lambda_ct} \\
    \mathbb{P}(Y \geqslant  t \,|\, F_2) &= \mathbb{P}(E_1 + E_2 > t) = (1+\lambda_bt)e^{-\lambda_bt}
\intertext{and}
    \mathbb{P}(X_a(0) = 2 \,|\, Y \geqslant  t, F_2) &= \mathbb{P}(E_1 > t) = e^{-\lambda_b t} \\
    \mathbb{P}(X_a(0) = 3 \,|\, Y \geqslant  t, F_2) &= \mathbb{P}(E_1 \leqslant t, E_1 + E_2 > t) =
    \int_0^t \mathbb{P}(E_2 \geqslant  t - s)\lambda_b e^{-\lambda_b s}ds = \lambda_b t e^{-\lambda_b t}.
\end{align*}
The probability that $X_a(0) = 4$ can be computed in a similar fashion.
}}\hfill$\Diamond$
\end{example}

\begin{example} \label{ex:recurrent}
{\em {
In our compact version of the {\sc mvm}, we solely
consider the impact of non-recurrent 
incidents on the vehicle speeds, and
ignore e.g.\ daily traffic patterns.
We are, however, able to capture these
recurrent events as well. To this end, 
we propose a similar strategy as in \cite{Levering2022Estimation}. 
There it is shown that the incident
duration, inter-incident time, and vehicle
speeds, are dependent on the time-of-day, 
but that these time-dependencies
can be tackled by working with periods of the day
$\Theta_1,\dots,\Theta_\ell$
over which these effects are essentially constant.
For every such period, the per-link incident
and inter-incident distribution are estimated,
as well as the corresponding driveable speed levels.
These can then be used when computing (and storing) the 
per-link travel-time distributions.

Representing recurrent events, the boundaries
of the periods $\Theta_1,\dots,\Theta_\ell$ 
are quite predictable.
Therefore, 
{the time between}
these boundaries
can be modeled by Erlang distributions.
Indeed, for given $k \in \mathbb{N}$ and
\mbox{$Z_i \sim \text{Exponential}( k/t)$,} $Z_1,\dots,Z_k$ independent, we have that $\sum_{i=1}^k Z_i$ is Erlang($k,k/t)$ distributed, and
\begin{align*}
\mathbb{E}\Big[\sum_{i=1}^k Z_i \Big] = t, \quad \quad \quad
\text{Var}\Big[\sum_{i=1}^k Z_i \Big] = t^2/k.
\end{align*}
Thus, modeling the time
between the boundaries of $\Theta_j$ with mean~$t_j$
by an Erlang-$k_j$ distribution with mean~$t_j$, 
we can achieve a low variance by choosing $k_j$ large enough.
% Now, to include $\Theta_1,\dots,\Theta_\ell$ in 
% the {\sc mvm}, we expand the background 
% process $B(t)$ with a Markov process $Y(t)$ that
% models the transitions of the Erlang phases 
% within and between  
% periods. 
Then, to include the periods $\Theta_1,\dots,\Theta_\ell$ in 
the {\sc mvm}, we expand the background 
process $B(t)$ with a Markov process $Y(t)$,
whose
state space consists of the
$k_1 + \dots + k_\ell$
Erlang phases that model
the times between their boundaries.
The process $Y(t)$ visits these states cyclically,
with $Y(t) = y$ encoding 
presence in period $\Theta_j$ at time $t$
if $y$ belongs to one of the $k_j$ Erlang
phases modeling the boundaries of $\Theta_j$.
Now, working with the extended $B(t)$,
we set the velocity of a vehicle
traversing $a_i \in A$ in the following way: the vehicle
speed at time $t$ equals $v_{a_i}(s_i,y)$ if
$X_{a_i}(t) = s_i \in \{1,2\}$ and
$Y(t) = y$. This way, the speed on a link depends
both on the time-of-day (via $Y(t)$) and 
the presence of an incident (via $X_{a_i}(t)$).

% Discuss: link and path distribution
As the velocity dynamics on $a_i$ are fully described
by the Markov process $(X_{a_i}(t),Y(t))$, we can
apply a discretization procedure in a similar
fashion as in Section~\ref{subsec:traveltimedist}
to obtain the link travel-time distribution on $a_i$.
Then, the list $L_{a_i}^{j,y}$ represents the travel
time distribution on $a_i$ given that $(X_{a_i}(t),Y(t))$
is initially in state $(j,y)$. 
We do, however, need to pay special attention to the
computation of the travel-time distribution
on paths, as the process
$Y(t)$ affects the velocities on \textit{all} arcs.
Therefore, instead of (travel time, probability)-pairs,
we let the lists $L_{a_i}^{j,y}$
contain (travel time, probability, state)-tuples,
in which the state is the state of $Y(t)$ in the final 
iteration step. Then, the travel time 
distribution on the path $\{a_1,\dots,a_{i+1}\}$, given that 
$L$ is the list of (travel time, probability, state)-tuples on
$\{a_1,\dots,a_{i}\}$, is simply given by the list
\begin{equation}
    \big\{\big(t_1 + t_2, p_1 \cdot p_2 \cdot [e^{t_1Q_{a_{i+1}}}]_{(s_{i+1}, j)}
    \big) \,|\, (t_1,p_1,y_1) \in L, 
    (t_2,p_2,y_2) \in L_{a_{i+1}}^{j,y_1}, j = 1,2\big\}.
    \label{eq:pathdistrecurrent}
\end{equation}

% initial state
The path-extension procedure in \eqref{eq:pathdistrecurrent} 
can now be used to compute the optimal departure time with 
one of the algorithms presented in Section~\ref{subsec:online}.
Note that, to use these algorithms, we need knowledge
on the state $Y(0)$. Similar to the previous phase-type
example, this state can not directly be observed, but
a probability distribution over the possible states of $Y(t)$ can be computed.
That is, given the request time, we do know
the current period~$\Theta_i$ and the elapsed time $t > 0$
between the start of this period and the request time.
Denoting with $y_1,\dots,y_k$ the subsequent Erlang states 
that model the duration of $\Theta_i$, and with $\lambda$ their transition
rate, we have
\begin{align*}
    \mathbb{P}(Y(0) = y_i \,|\, Y(-t) = y_1, Y(0) \in \{y_1,\dots,y_k\}) &=
    \frac{\mathbb{P}(Y(0) = y_i \,|\, Y(-t) = y_1)}
    {\mathbb{P}(Y(0) \in \{y_1,\dots,y_k\} \,|\, Y(-t) = y_1)} \\
    &= \frac{\mathbb{P}(S_{i-1} < t, S_i \geqslant  t)}{\mathbb{P}(S_k \geqslant  t)},
\end{align*}
with $S_j = \text{Erlang}(j, \lambda)$. Conditioning on the
value of $S_{i-1}$, it is now easily derived that, 
for $E_i \sim \text{Exponential}(\lambda)$,
\begin{align*}
    \frac{\mathbb{P}(S_{i-1} < t, S_i \geqslant  t)}{\mathbb{P}(S_k \geqslant  t)} &=
    \frac{\int_{0}^t \mathbb{P}(E_i > t\!-\!s)f_{S_{i-1}}(s)ds}{\mathbb{P}(S_k \geqslant  t)}
    = \frac{e^{-\lambda t} \int_{0}^t \frac{(\lambda s)^{i-2}}{(i-2)!} ds}{\mathbb{P}(S_k \geqslant  t)}
    = \frac{(\lambda t)^{i-1}}{\lambda (i\!-\!1)! \sum_{n = 0}^{k-1} \frac{1}{n!}(\lambda t)^n}.
\end{align*}
}}\hfill$\Diamond$
\end{example}

%% Alternative method
% One drawback of incorporating recurrent events via
% Erlang phases is the increase in complexity, as the
% state space size of $B(t)$ is relatively large.
% Therefore, we also describe a heuristic alternative, 
% which simply stores (travel time, probability)-pairs
% per time slot $\Theta_i$, and uses a bisection
% algorithm per time slot to deduce the optimal departure
% time. The general idea of this method is described
% in Appendix~\ref{appendix:alternativerecurrent}. \textit{Include
% this or not?}

\section{Optimal departure time: speed-up techniques}\label{app:speedup}~\\
We will discuss a simple but effective way
to speed-up the computation of the optimal departure
time for certain paths or OD pairs. 
% The idea is to not only precompute per-link
% travel-time distributions, but to precompute
% the travel-time distributions of a set of (frequently used)  paths as well.
% Importantly, for long
% path lengths or OD-pair distances, the computational
% savings of these additional precomputations
% are typically substantial, and will thus allow
% the practical use of the algorithms in larger networks.
% present a way to reduce the run-time 
% of Algorithm 2 and the k-shortest path
% method in case the paths consist of a high number of links
% Idea: pre-compute for often traveled OD-pairs ('hubs/landmarks')
We have limited the computational costs by storing
per-link travel-time distributions.
%which, for a given link, are lists of (travel time, probability)-pairs 
%for the initial states the Markov process on the link can attain.
We could, however, expand the work that is carried out in
the precomputations and store the travel-time distribution
for some sets of subsequent links as well. 
Concretely, for a given path~$P'$, we are able to compute
the list $L_{P'}^j$ of (travel time, probability)-pairs that arises
for each initial state $j$ the background process of $P'$ can
attain. Then, when computing the optimal departure time
on a path~$P$ that has $P'$ as subpath, we can simply
view $P'$ as one link and use the stored travel
time distributions for $P'$ in the algorithms. 
Notably, this will speed up the computation of the
optimal departure time in both Algorithm~\ref{alg:offline}
and the $k$-shortest path method,
whenever (one of) the path(s) contains subpath(s) for which
the travel-time distribution is stored.
In contrast, as the A$^\star$-algorithm works with
individual links rather than with paths, the speed-up
technique is not applicable to Algorithm~\ref{alg:od1}.

\begin{remark}
{\em 
Availability of the precomputed travel-time distribution of a subpath $P'$
reduces the number of iterations necessary to compute the travel-time distribution
on a complete path~$P$. Note that this does not directly yield a speed-up, 
as the number of operations within the iteration adding $P'$ is
potentially large. That is, with $P'$ consisting of $k$ individual
links, this iterative step works -- instead of with two -- with 
$2^k$ (travel time, probability) lists, 
corresponding to the number of states in
the background process of $P'$.
% this iterative step is similar
% to the step presented in
% \eqref{eq:pathdist}, but, corresponding to the number of states in
% the background process of $P'$, allows $j \in \{1,2,\dots,2^k\}$.
However, as the operations
concerning different background states can be carried out
in parallel, the impact of these additional operations is greatly
reduced.  
% The fact that the additional precomputations \textit{do} indeed reduce
% the run-time, stems from the observation that the operations
% concerning different background states can be carried out
% in parallel.
}\hfill$\Diamond$
\end{remark}

% Reduce state space by deleting small probabilities
There is an additional advantage when precomputing
the travel-time distribution for a set of paths.
That is, with only a negligible loss in accuracy,
we are able to substantially reduce the size of the
state space and, consequently, further speed
up the computations. Concretely, if one wants
to find the optimal departure time
for a path $P$ containing a subpath $P'$
for which the distribution is stored, and if there are states 
in the state space of 
%the joint Markov process on 
$P'$
that have an extremely small probability of occurring, then the distributions
that correspond to these states can simply be neglected. 
This means that, in the
computation of the travel-time distribution on the
path~$P$, we will omit lists $L_{P'}^j$
of (travel time, probability)-pairs on $P'$
for which 
\begin{align*}
    \mathbb{P}(\exists t \in [0,M]: (X_{a'}(t))_{a' \in P'} = j \,|\,
    B(0) = s) < \varepsilon,
\end{align*}
for some small $\varepsilon > 0$.
Examples of procedures to generate upper bounds for such hitting probabilities
are provided in \cite{Levering2021AConditions}.

\end{document}